\documentclass[11pt]{amsart}
\usepackage{amsfonts,amssymb}
\usepackage[all]{xy}

\numberwithin{equation}{section}

\newtheorem{theo}{Theorem}[section]
\newtheorem{prop}[theo]{Proposition}
\newtheorem{lem}[theo]{Lemma}
\newtheorem{coro}[theo]{Corollary}

\theoremstyle{definition}
\newtheorem{defi}[theo]{Definition}
\newtheorem{example}[theo]{Example}

\newtheorem{rem}[theo]{Remark}

\newtheorem{notation}[theo]{Notation}
\newtheorem{convention}[theo]{Convention}

\newcommand\N{{\mathbb{N}}}

\newcommand\T{{\mathbb{T}}}
\newcommand\aA{{\mathcal{A}}}
\newcommand\bB{{\mathcal{B}}}
\newcommand\cC{{\mathcal{C}}}

\newcommand\jJ{{\mathcal{J}}}
\newcommand\lL{{\mathcal{L}}}
\newcommand\nN{{\mathcal{N}}}
\newcommand\pP{{\mathcal{P}}}
\newcommand\qQ{{\mathcal{Q}}}
\newcommand\uU{{\mathcal{U}}}
\newcommand\vV{{\mathcal{V}}}
\newcommand\wW{{\mathcal{W}}}
\newcommand\aAb{{\mathcal{A}}^{\scriptscriptstyle{\bullet}}}
\newcommand\bBb{{\mathcal{B}}^{\scriptscriptstyle{\bullet}}}

\newcommand\uUb{{\mathcal{U}}_{\scriptscriptstyle{\bullet}}}
\newcommand\vVb{{\mathcal{V}}_{\scriptscriptstyle{\bullet}}}
\newcommand\wWb{{\mathcal{W}}_{\scriptscriptstyle{\bullet}}}
\newcommand\Mb{M_{\scriptscriptstyle{\bullet}}}
\newcommand\eps{\varepsilon}
\newcommand\teps{\tilde{\varepsilon}}
\newcommand\teta{\tilde{\eta}}
\newcommand\tf{\tilde{f}}
\newcommand\tg{\tilde{g}}
\newcommand\tp{\tilde{p}}
\newcommand\vH{\check{H}}
\newcommand\vC{\check{C}}

\newcommand\id{{\mathrm{Id}}\,}
\renewcommand\hom{{\mathrm{Hom}}\,}
\newcommand\et{\'{e}tale}
\newcommand\inferieur{<}

\newcommand\com{{\scriptscriptstyle{\bullet}}}
\newcommand\equal{=}

\newcommand\pf{\begin{proof}}
\newcommand\pfend{\end{proof}}

\begin{document}
\title{Groupoid cohomology and extensions}
\author{by J.~L.~Tu}

\date{}

\begin{abstract}
We show that Haefliger's cohomology for \'etale groupoids, Moore's
cohomology for locally compact groups and the Brauer group of
a locally compact groupoid are all particular cases of sheaf (or
\v{C}ech) cohomology for topological simplicial spaces.
\end{abstract}
\maketitle
\tableofcontents
\section{Introduction}
Let $G$ be a locally compact Hausdorff
groupoid with Haar system. In \cite{KMRW98},
the authors studied the group of Morita equivalence classes
of actions of $G$ on continuous fields of $C^*$-algebras over
the unit space $G_0$ such that
\begin{itemize}
\item Each fiber is isomorphic to the algebra of compact operators
on some Hilbert space (depending on the fiber);
\item The bundle satisfies Fell's condition, i.e.
each point of $G_0$ has a neighborhood $U$ such that there exists
a section $f(x)$ with $f(x)$ a rank-one projection for all
$x\in U$.
\end{itemize}
They called this group the Brauer group $Br(G)$ of $G$, and showed that
it is naturally isomorphic to the group $\mathrm{Ext}(G,\T)$
of Morita equivalence classes of central extensions
$$\T\times G'_0 \to E\to G',$$
where $G'$ is some Morita equivalent groupoid. In the case of discrete
groups, it is well-known that central extensions of $G$ by $\T$
are classified by $H^2(G,\T)$. Actually, given any locally compact
group $G$ and any Polish (i.e. metrizable separable complete)
$G$-module $A$, Moore's cohomology groups $H^2(G,A)$
classify extensions $A\hookrightarrow E\twoheadrightarrow G$
such that the action of $G$ on $A$ by conjugation is exactly
the action of $G$ on the $G$-module $A$ \cite{Moo64,Moo76}.
One of the possible definitions of Moore's cohomology is the
following: consider $C^n(G,A)$ the space of all measurable maps
$c\colon G^n\to A$ with the differential
\begin{eqnarray}\label{eqn:moore-differential}
\lefteqn{(dc)(g_1,\ldots,g_n){\equal}g_1c(g_2,\ldots,g_n)}\\\nonumber
&&+\sum_{k{\equal}1}^n (-1)^k
c(g_1,\ldots,g_kg_{k+1},\ldots,g_n)+(-1)^{n+1}c(g_1,\ldots,g_n),
\end{eqnarray}
then $H^n(G,A)$ is the $n$-th cohomology group of the complex
$C^*(G,A)$.
\par\medskip

On the other hand, Haefliger (\cite{Hae79}, see also
\cite{Kum88}) defined sheaf cohomology groups
$H^*(G,\aA)$ given any \et\ groupoid $G$ and any abelian $G$-sheaf $\aA$
(i.e. an abelian sheaf on $G_0$ endowed with a continuous action of $G$).
It was thus natural to expect that a single cohomology theory for groupoids
should unify all these: this is the question asked by A.~Kumjian in
Boulder (1999).
\par\bigskip

Our approach is to consider the simplicial (topological) space
$G_\com{\equal}(G_n)_{n\in \N}$ and to use sheaf cohomology for simplicial
spaces \cite{Del74} and \v{C}ech cohomology (see Section~\ref{sec:cech}).
We show:

\begin{theo}
Let $\Mb$ be a simplicial space and $\aAb$ be an abelian sheaf
on $\Mb$. Denote by $H^n(\Mb;\aAb)$ and $\vH(\Mb;\aAb)$ the
sheaf and \v{C}ech cohomology groups respectively. Then
\begin{itemize}
\item[(a)] $H^n(\Mb;\aAb)\cong \vH(\Mb;\aAb)$ if $M_n$ is paracompact
for all $n$;
\item[(b)] $H^n(G_\com;\aAb)$ coincides with Haefliger's cohomology
$H^n(G;\aA^0)$ if $G$ is an \et\ groupoid and $\aAb$ is the sheaf on
$G_\com$ corresponding to an abelian $G$-sheaf $\aA^0$;
\item[(c)] $H^n(G_\com;\aAb)\cong \vH(G_\com;\aAb)
\cong H^n_{\mathrm{Moore}}(G,A)$
if $G$ is a locally compact group, $A$ is a Polish $G$-module
and $\aAb$ is the sheaf on $G_\com$ associated to $A$;
\item[(d)] $H^2(G_\com;\T)\cong \vH(G_\com;\T)$
is the Brauer group of $G$ if $G$ is a locally
compact Hausdorff groupoid with Haar system;
\item[(e)] $H^n(G_\com,\aAb)$ and $\vH^*(G_\com,\aAb)$ are invariant
under Morita equivalence of topological groupoids.
\end{itemize}
\end{theo}

(Actually we prove more than statement (d);
see Proposition~\ref{prop:2-cohomology}.)

It is possible that some parts
of this paper are well-known among specialists,
but apparently written nowhere in the literature, in particular
the definition of \v{C}ech cohomology and its relation to
sheaf cohomology for simplicial spaces. Besides, a more conceptual
approach would have consisted in using Grothendieck's cohomology
of sites \cite{Art72}, as Moerdijk did in \cite{Moe91}.
Moreover, it is possible that the present work has non-trivial
intersection with Moerdijk's in \cite{Moe03}.
However, we hope that the present
approach, being rather direct and elementary,
will still be of interest to the reader.

{\bf Acknowledgments:} the author would like to thank Ping Xu for
useful discussions and Kai Behrend for providing some bibliographic
references.

\section{Simplicial spaces and groupoids}
\subsection{Definition of simplicial spaces}
Let us recall some basic facts about simplicial spaces.
Let $\Delta$ (resp. $\Delta'$) be the simplicial (resp. semi-simplicial)
category, whose objects are the nonnegative integers, and whose
morphisms are the nondecreasing (resp. increasing) maps
$[m]\to [n]$ (where $[n]$ denotes the interval $\{0,\ldots,n\}$.
We denote by $\Delta^{(N)}$ the $N$-truncated simplicial
category, i.e. the full sub-category of $\Delta$ whose objects
are the integers $\le N$.
\par\medskip

A simplicial (resp. semi-simplicial, $N$-simplicial) topological space
is a contravariant functor from the category
$\Delta$ (resp. $\Delta'$, resp. $\Delta^{(N)}$)
to the category of topological spaces. In the same way,
one can define the notion of simplicial (resp. semi-simplicial, $N$-simplicial)
manifold. In this paper we shall work with simplicial topological
spaces and will use the terminology ``simplicial space'', but
some of the results can easily be transposed to simplicial
manifolds.
\par\medskip

In practice, a (semi-)simplicial space is a sequence
$M_\com{\equal}(M_n)_{n\in\N}$ of topological spaces, given with
continuous maps $\tilde f\colon M_n\to M_k$ for every morphism
$f\colon [k]\to [n]$, satisfying the relation
$\widetilde{f\circ g}{\equal}\tilde g\circ \tilde f$ for all composable
morphisms $f$ and $g$.
\par\medskip

Let $\eps^n_i\colon [n-1]\to [n]$ be the unique increasing map that
avoids $i$, and $\eta^n_i\colon [n+1]\to [n]$ be the unique non-decreasing
surjective map such that $i$ is reached twice ($0\le i\le n$).
We will usually omit the superscripts for convenience of notation.

If $M_\com{\equal}(M_n)_{n\in\N}$ is a simplicial space, then the
face maps $ \teps_i^n: M_n \to M_{n-1}$, $i{\equal} 0,\dots,n$
and the degeneracy maps $ \teta_i^n: M_n \to M_{n+1}$, $i{\equal} 0,\dots,n$,
satisfy the following simplicial identities:
$\teps_i^{n-1} \teps_j^{n}{\equal}  \teps_{j-1}^{n-1} \teps_i^{n}$
if $i<j$,
$\teta_i^{n+1} \teta_j^{n}{\equal} \teta_{j+1}^{n+1} \teta_i^{n} $ if $i \leq j$,
$\teps_i^{n+1} \teta_j^{n}{\equal}  \teta_{j-1}^{n-1} \teps_{i}^{n} $
if $i < j $, $\teps_i^{n+1} \teta_j^{n}{\equal}
\teta_{j}^{n-1} \teps_{i-1}^{n} $
if $i > j+1 $ and $\teps_j^{n+1} \teta_j^{n}{\equal}
\teps_{j+1}^{n+1} \teta_j^{n}{\equal} \id_{M_n} $.
\par\medskip
Conversely, if we are given a sequence $M_\com$ of topological
spaces and maps satisfying such identities, then there is a unique
simplicial structure on $M_\com$ such that $\teps_i^n$ are the
face maps and $\teta^n_i$ are the degeneracy maps.

\subsection{Groupoids}
In order to fix notations,
we first recall some basic facts about groupoids. For more details,
see e.g. \cite{Ren80}.

A topological groupoid is given by two topological spaces
$G_0$ and $G$, two maps $r$ and $s$ from $G$ to $G_0$, called
the range and source maps, a unit map $\eta\colon G_0\to G$,
a partially defined multiplication $G_2{\equal}\{(g,h)\in G^2|\;
s(g){\equal}t(h)\} \to G$ denoted by $(g,h)\mapsto gh$, and an inversion
map $G\to G$ denoted by $g\mapsto g^{-1}$ such that the following
identities hold (for $g,h,k\in G$ and $x\in G_0$):
\begin{itemize}
\item $r(gh){\equal}r(g)$, $s(gh){\equal}s(h)$;
\item $(gh)k{\equal}g(hk)$ whenever $s(g){\equal}r(h)$ and $s(h){\equal}r(k)$;
\item $s(\eta(x)){\equal}r(\eta(x)){\equal}x$;
\item $g\eta(s(g)){\equal}\eta(r(g)) g {\equal} g$;
\item $r(g^{-1}){\equal}s(g)$, $s(g^{-1}){\equal}r(g)$, $gg^{-1}{\equal}\eta(r(g))$,
$g^{-1}g{\equal}\eta(s(g))$.
\end{itemize}

We will usually identify the unit space $G_0$ to a subspace of $G$
by means of the unit map $\eta$.

Standard examples are:
\begin{itemize}
\item groups, with $G_0{\equal}\mathrm{pt}$;
\item spaces $M$, with $G{\equal}G_0{\equal}M$, $r{\equal}s{\equal}\mathrm{Id}$;
\item the homotopy groupoid of a space $M$, where $G_0{\equal}M$,
$G$ is the set of homotopy classes of paths in $M$, $s(g)$ is the
starting point of the path and $r(g)$ is the endpoint.
\end{itemize}

Here are a few notations that we will use: $G_x{\equal}s^{-1}(x)$,
$G^x{\equal}r^{-1}(x)$, $G_x^y{\equal}G_x\cap G^y$.
\par\medskip

A left \emph{action} of a topological groupoid $G$ on a space $Z$
is given by a
(continuous) map $p\colon Z\to G_0$ and a map
$G\times_{s,p}Z\to Z$, denoted by $(g,z)\mapsto gz$, such that
\begin{itemize}
\item $p(gz){\equal}r(g)$;
\item $(gh)z{\equal}g(hz)$ whenever $(g,h)\in G_2$ and $s(h){\equal}p(z)$;
\item $ez{\equal}z$ if $e\in G_0\subset G$.
\end{itemize}

We will say that $Z$ is a (left) $G$-space.
Given a $G$-space $Z$, we form the crossed-product
groupoid $G\ltimes Z:{\equal}G\times_{s,p} Z$ with unit space $Z$,
source and range maps $s(g,z){\equal}z$, $r(g,z){\equal}gz$, product $(g,z)(h,z')
{\equal}(gh,z')$ if $z{\equal}hz'$, and inverse $(g,z)^{-1}{\equal}(g^{-1},gz)$.
\par\bigskip

Any topological groupoid $G$ canonically gives rise to a simplicial space
as follows \cite{Seg68}:
Let
$$G_n{\equal}\{(g_1,\ldots,g_n)\vert\;
s(g_i){\equal}t(g_{i+1})\;\forall i\}$$
be the set composable $n$-tuples.

Define the face maps $\teps_{i}^n :G_{n}\to
G_{n-1}$ by, for $n > 1$
 \begin{eqnarray*} 
 && \teps_{0}^n (g_1  ,
 g_2 , \ldots , g_n){\equal}  ( g_2 , \ldots ,g_n)\\
  && \teps_{n}^n (g_1  , g_2 ,
 \ldots , g_n){\equal}  (g_1 , \ldots ,g_{n-1})\\
 && \teps_{i}^n (g_1  ,  \ldots,  g_n){\equal} 
 (g_1 ,  \ldots, g_i g_{i+1} , \ldots ,g_n), \ \  1\leq
i \leq n-1,
\end{eqnarray*}
and for $n{\equal} 1$ by,
$ \teps_{0}^1 (g){\equal} s (g)$,
$\teps_{1}^1 (g){\equal} r (g)$.
Also define the degeneracy maps: $\teta_0^0
: G_0 \to G_1$ is the unit map of the groupoid, and
$\teta^n_i :G_{n}\to G_{n+1}$ by:
\begin{eqnarray*}
&& \teta^n_0 (g_1 ,\ldots , g_n){\equal} 
(r ( g_1), g_1 , \ldots , g_n) \\
&& \teta^n_i (g_1 ,\ldots , g_n){\equal} 
(g_1 , \ldots ,g_{i}, s (g_{i})
,g_{i+1}, \ldots , g_n), \ \ 1\leq i\leq n.
\end{eqnarray*}

Another way to view the simplicial structure of $G_\com$
is the following: we note that $G_n$ can be identified
with the quotient of
$(EG)_n:{\equal}\{(\gamma_0,\ldots,\gamma_n)\in G^{n+1}\vert\;
r(\gamma_0){\equal}\cdots{\equal}r(\gamma_n)\}$ by the left action of $G$,
the correspondence being
\begin{eqnarray*}
(g_1,\ldots,g_n)
&{\equal}&(\gamma_0^{-1}\gamma_1,\ldots,\gamma_{n-1}^{-1}\gamma_n)\\
{[} \gamma_0,\ldots,\gamma_n ] &{\equal}& [ r(g_1),g_1,g_1g_2,\ldots,g_1\ldots g_n ].
\end{eqnarray*}
Then, for any morphism $f\colon [k]\to [n]$, $\tf\colon G_n\to G_k$
is defined by
$$\tf[\gamma_0,\ldots,\gamma_n]
{\equal}[\gamma_{f(0)},\ldots,\gamma_{f(n)}].$$
For instance, in the first picture,
if $f$ is injective then
\begin{equation}\label{eqn:groupoid-tilde-f}
\tf(g_1,\ldots,g_n){\equal}
(g_{f(0)+1}\cdots g_{f(1)}, g_{f(1)+1}\cdots g_{f(2)},\ldots,
g_{f(k-1)+1}\cdots g_{f(k)}).
\end{equation}

\subsection{Morita equivalence and generalized morphisms}\label{sec:Morita}
We recall (see for instance \cite{Hae84,Hil-Sk87,Leg99,Mrc99,Lan01,TXL03})
that a generalized morphism between two (topological,
or locally compact, or Lie) groupoids
$G'$ and $G$ is given by a topological space (or a locally compact
space, or a manifold) $Z$, two maps
$G'_0\stackrel{\rho}{\leftarrow}
Z\stackrel{\sigma}{\to}G_0$ such that $Z$ admits a left action of
$G'$ with respect to $\rho$, a right action of $G$ with respect to
$\sigma$, with the property that the two actions commute and
$\rho\colon Z\to G'_0$
is a locally trivial $G$-principal bundle.

Topological (or locally compact...)
groupoids and generalized morphisms form a category
whose isomorphisms are Morita equivalences. Every groupoid morphism
naturally defines a generalized morphism.

If $\uU{\equal}(U_i)_{i\in I}$ is an open cover of $G_0$, define the cover groupoid
\begin{equation}\label{eqn:GU}
G[\uU]{\equal}\{(i,g,j)\in I\times G\times I
\vert\; r(g)\in U_i, s(g)\in U_j\}
\end{equation}
with unit space $\{(i,x)\in I\times G_0\vert\; x\in U_i\}$,
source and range maps $s(i,g,j){\equal}(j,s(g))$, $r(i,g,j){\equal}(i,r(g))$
and product $(i,g,j)(j,h,k){\equal}(i,gh,k)$.

Then the canonical morphism $G[\uU]\to G$ is a Morita equivalence.
Moreover, every generalized morphism $G'\to G$ admits a
decomposition $G'\stackrel{\sim}{\leftarrow}G'[\uU']
\stackrel{f}{\to} G$ for some open cover $\uU'$ of $G'_0$ and some
\emph{groupoid morphism} $f$.

\begin{rem}\label{rem:GU}
The simplicial space $G[\uU]_\com$ is isomorphic to the sub-simplicial
space of $(I^{n+1}\times G_n)_{n\in \N}$ such that $G[\uU]_n$
consists of $(2n+1)$-tuples $(i_0,\ldots,i_n,g_1,\ldots,g_n)$
satisfying the condition
\begin{equation}\label{eqn:GUn}
r(g_1)\in U_{i_0},\;s(g_1)\in U_{i_1},\ldots,s(g_n)\in U_{i_n}.
\end{equation}
\end{rem}

Finally, we will need the following
\begin{prop}\label{prop:Morita}
For any functor $F$ from the category of topological
groupoids to any category,
the following are equivalent:
\begin{itemize}
\item[(i)] $F$ is invariant under Morita-equivalence;
\item[(ii)] $F$ factors through the category whose objects
are groupoids and whose morphisms are generalized morphisms;
\item[(iii)] For any groupoid $G$ and any open cover $\uU$ of $G_0$,
the canonical map $G[\uU]\to G$ induces an isomorphism
$F(G[\uU])\stackrel{\sim}{\to} F(G)$.
\end{itemize}
\end{prop}

\pf
See for instance \cite[Proposition~2.5]{TXL03}.
\pfend

\section{Sheaves on simplicial spaces}
\subsection{Basic definitions}
Recall \cite{Del74} that if $u\colon X\to Y$ is continuous,
$\aA$ is a sheaf on $X$ and $\bB$ is a sheaf on $Y$, then
a $u$-morphism from $\bB$ to $\aA$ is by definition an element
of $\hom(\bB,u_*\aA)\cong \hom(u^*\bB,\aA)$. A sheaf on a
simplicial (resp. semi-simplicial) space $M_\com$ is a sequence
$\aA^\com{\equal}(\aA^n)_{n\in\N}$ such that $\aA^n$ is a sheaf on $M_n$,
and such that for each morphism $f\colon [k]\to [n]$ in the
category $\Delta$ (resp. $\Delta'$) we are given $\tf$-morphisms
\begin{equation}\label{eqn:f*}
\tf^*\colon \aA^k\to \aA^n
\end{equation}
such that $\tf^*\tg^*
{\equal}\widetilde{f\circ g}^*$ if $g\colon [\ell]\to [k]$

In practice, given open sets $U\subset M_n$ and $V\subset M_k$
such that $\tf(U)\subset V$ we have a \emph{restriction map}
$\tf^*\colon \aA^k(V)\to \aA^n(U)$ such that
$\tf^*\circ\tg^*{\equal}\widetilde{f\circ g}^*\colon
\aA^\ell(W)\to\aA^n(U)$ whenever $\tg(V)\subset W$.

A fundamental example is given by $G$-sheaves. In the definition below,
recall that a map $f\colon X\to Y$ is said to be \et\ if it is
a local homeomorphism, i.e. every point $x\in X$ has an open neighborhood
$U$ such that $f(U)$ is open and $f$ induces a homeomorphism
from $U$ onto $f(U)$. We will also say that $X$ is an \et\ space over $Y$.
A groupoid is \et\ if the range (equivalently the source) map is
\et. A morphism $\pi_\com\colon E_\com\to M_\com$
is \et\ if each $\pi_n\colon E_\com\to M_\com$ is \et.
Finally, recall that a sheaf over a space $X$ can be considered as
a (not necessarily Hausdorff) \et\ space over $X$.

\begin{defi}\cite{Hae79}
Let $G$ be a topological groupoid. Then a \emph{$G$-sheaf} is an \et\ space
$E_0$ over $G_0$, endowed with a continuous action of $G$.
\end{defi}

Of course, an \emph{abelian} $G$-sheaf is a $G$-sheaf $E_0$ such that
$E_0$ is an abelian sheaf on $G_0$ and such that for each
$g\in G$, the action $\alpha_g\colon (E_0)_{s(g)}\to (E_0)_{r(g)}$
is a group morphism.

\begin{example}
If $G$ is a group then a $G$-sheaf is just a space endowed with
an action of $G$.
\end{example}

To show that any $G$-sheaf defines a sheaf over $G_\com$,
we need some preliminaries:

\begin{defi}\label{defi:reduced}
Let $\pi_\com\colon E_\com\to M_\com$ a morphism
of simplicial spaces. We say that $\pi_\com$ is \emph{reduced} if
for all $k,n$ and all $f\in\hom_\Delta(k,n)$,
the map $\tf$ induces an isomorphism
$E_n\cong M_n\times_{\tf,\pi_k}E_k$.
In this case, we will say that $E_\com$ is a reduced simplicial
space over $M_\com$.
\end{defi}

\begin{defi}\label{defi:reduced2}
Let $\aA^\com$ be a sheaf over the simplicial space $M_\com$.
We will say that $\aA^\com$ is \emph{reduced} if
for all $k,n$ and all $f\in\hom_\Delta(k,n)$,
the morphism $\tf^*\in
\hom(\tf^*\aA^k,\aA^n)$ is an isomorphism.
\end{defi}

\begin{lem}\label{lem:reduced}
There is a one-to-one correspondence between
reduced sheaves over $M_\com$ and reduced \et\ simplicial
spaces over $M_\com$.
\end{lem}

\pf
The proof is easy. Let us just explain the construction of the sheaf
$\aA^\com$ out of the reduced simplicial space $E_\com$
over $M_\com$.

Let $\aA^n(U)$ be the space of
continuous sections over $U$ of the projection map $\pi_n\colon E_n\to M_n$.
If $f\colon [k]\to [n]$ is a morphism in $\Delta$
and $\tf(U)\subset V$, then
for any section $\sigma\in \aA^k(V)$ we define
$\tf^*\sigma\in \aA^n(U)$ by
$(\tf^*\sigma)(x){\equal}(x,\sigma(\tf(x)))
\in M_n\times_{\tf,\pi_k}E_k\cong E_n$.
\pfend

\begin{lem}\label{rem:reduced2}
Any reduced simplicial space over $M_\com$,
\et\ or not, determines a sheaf over $M_\com$.
\end{lem}

\pf
The proof is the same. Note that it is not clear
whether all sheaves can be constructed this way.
\pfend

\begin{coro}\label{coro:reduced}
Let $G$ be a topological groupoid, then
any $G$-space determines a sheaf on $G_\com$.
If the $G$-space is \et\ then it determines a reduced sheaf on
$G_\com$.
\end{coro}

\pf
Suppose $Z$ is a $G$-space and let $\pi_n$ be the first projection
$(G\ltimes Z)_n{\equal}G_n\times_{\tp_n,p}Z\to G_n$,
where $\tp_n(g_1,\ldots,g_n){\equal}s(g_n)$. Then $\pi_n$ is clearly a
simplicial map $(G\ltimes Z)_\com \to G_\com$.
\pfend

\begin{coro}\label{coro:G-sheaf-simplicial}
Every $G$-sheaf canonically defines a reduced sheaf over the simplicial
space $G_\com$.
\end{coro}

Another example is given by $G$-modules:

\begin{defi}\label{def:G-module}
Let $G$ be a topological groupoid. A $G$-module is a topological
groupoid $A$,
with source and range maps equal to a map $p\colon G_0$, such that
\begin{itemize}
\item $A_x^x$ is an abelian group for all $x$;
\item As a space, $A$ is endowed with a $G$-action $G\times_{s,p}A\to A$;
\item for each $g\in G$, the map $\alpha_g\colon A_{s(g)}\to A_{r(g)}$
given by the action is a group morphism.
\end{itemize}
\end{defi}

By Corollary~\ref{coro:reduced}, any $G$-module
defines a sheaf $\aA^\com$ which is clearly abelian.

More explicitly, the simplicial structure on $(G\ltimes A)_\com$
is defined as follows: for all $f\in\hom_\Delta(k,n)$,
\begin{equation}\label{eqn:simplicial-GxA}
\tf([\gamma_0,\ldots,\gamma_n],a){\equal}
([\gamma_{f(0)},\ldots,\gamma_{f(k)}],\gamma_{f(k)}^{-1}\gamma_n a).
\end{equation}

Then $\aA^n$ is the sheaf of germs of continuous sections
of $(G\ltimes A)_n \to G_n$, i.e. sections are continuous maps
$\varphi(g_1,\ldots,g_n)\in A_{s(g_n)}$. However, to recover the
usual formulas like (\ref{eqn:moore-differential}), it is better
to work with the maps
$$c(g_1,\ldots,g_n){\equal}g_1\cdots g_n\varphi(g_1,\ldots,g_n)\in
A_{r(g_1)}$$
and this is what we shall usually do.

Note that for all $\vec{g}{\equal}(g_1,\ldots,g_n)\in G_n$, the stalk
$\aA_{\vec{g}}$ maps to $A_{\tp_n(\vec{g})}$. This map is surjective
iff $p\colon A\to G_0$ has enough cross-sections; for injectivity,
it is enough that $p$ be an \et\ map.

If $A{\equal}G_0\times B$ has constant fibers ($B$ being a topological
abelian group with no action of $G$),
then the corresponding sheaf is called the
constant sheaf and is again (abusively) denoted by $B$.

When $G$ is a group, a $G$-module is just a topological abelian group $A$
endowed with a continuous action $G\to {\text{Aut}}(A)$ and
the sheaf $\aA^n$ is just the constant sheaf $A$ on $G_n$.

\subsection{$G$-sheaves and sheaves over simplicial spaces}
Recall (Corollary~\ref{coro:G-sheaf-simplicial}) that any
$G$-sheaf defines a sheaf over $G_\com$.
In this subsection, which can safely be omitted by the reader,
we examine the converse. More precisely:

\begin{prop}
Let $G$ be a topological groupoid.
There is a one-to-one correspondence between:
\begin{itemize}
\item[(i)] $G$-sheaves;
\item[(ii)] reduced sheaves over $G_\com$;
\item[(iii)] reduced \et\ spaces over $G_\com$.
\end{itemize}
\end{prop}

\pf
(ii) $\iff$ (iii) was proved in
Lemma~\ref{lem:reduced} and (i) $\implies$ (ii) is
given by Corollary~\ref{coro:G-sheaf-simplicial}.
For (iii) $\implies$ (i), the proposition below will allow us to conclude.
\pfend

\begin{prop}
Let $G$ be a topological groupoid. There is a one-to-one
correspondence between reduced morphisms of simplicial spaces
$\pi_\com\colon E_\com\to G_\com$ and $G$-spaces.
Under this correspondence, \et\ spaces over $G_\com$ are
mapped onto $G$-spaces which are \et\ over $G_0$.
\end{prop}

\pf
The only difficulty is to show that $E_0$ is endowed with an action
of $G$ such that $(G\ltimes E_0)_\com\cong E_\com$.
Consider the map
$$\phi_n\colon E_n\to G_n\times E_0^{n+1}$$
defined by $\phi_n{\equal}(\pi_n,\tp_0,\ldots,\tp_n)$
where $p_i\colon [0]\to [n]$ is the map $p_i(0){\equal}i$.

One can view $(E_0^{n+1})_{n\in\N}$ as a simplicial space,
with
$$\tf(\xi_0,\ldots,\xi_n){\equal}(\xi_{f(0)},\ldots,\xi_{f(k)})\quad
\forall f\in\hom_\Delta(k,n).$$

It is not hard to show,
using the relation $\tp_j\circ \tf {\equal} \tp_{f(j)}$ (which holds
because $f\circ p_j{\equal}p_{f(j)}$), that $\phi_\com{\equal}(\phi_n)_{n\in\N}$
is a morphism of simplicial spaces. 

Since $E_\com$ is reduced, $(\pi_n,\tp_n)$ is a homeomorphism
from $E_n$ to $G_n\times_{\tp_n,\pi_0}E_0$ and in particular
$\phi_n$ is injective.

Define the action of $G$ on $E_0$ be the composition of maps
$$G\times_{s,p}E_0\stackrel{(\pi_1,\tp_1)^{-1}}{\longrightarrow}
E_1\stackrel{\tp_0}{\to} E_0,$$
i.e.
\begin{equation}\label{eqn:action}
\xi'{\equal}g\xi \iff \exists x\in E_1,\;
\phi_1(x){\equal}(g,\xi',\xi).
\end{equation}

Let $x\in E_n$, let $(g_1,\ldots,g_n){\equal}\pi_n(x)$ and $\xi{\equal}\tp_n(x)$.
Let us show that
\begin{equation}\label{eqn:phi-n}
\phi_n(x){\equal}(\pi_n(x),g_1\cdots g_n\xi,g_2\cdots g_n\xi,\ldots,\xi),
\end{equation}

where $g_1\cdots g_n\xi {\equal} g_1(g_2(\cdots(g_n \xi)\cdots))$ (we don't know
yet that $(g,\xi)\mapsto g\xi$ is an action).

Consider the map $f\colon [1]\to [n]$ such that $f(0){\equal}j-1$ and $f(1){\equal}j$.
Let $y{\equal}\tf(x)\in E_1$, then
$$\phi_1(y){\equal}\phi_1(\tf(x)){\equal}\tf(\phi_n(x))
{\equal}(g_j,\tp_{j-1}(x),\tp_j(x)),$$
hence by (\ref{eqn:action}),
$\tp_{j-1}(x){\equal}g_j\tp_j(x)$. Equation~(\ref{eqn:phi-n}) follows
by reverse induction on $j$.
\par\medskip

Now, let $(g,h,\xi)\in G_2\times_{\tp_2,p}E_0$. We want to show that
$g(h\xi){\equal}(gh)\xi$. Let $x{\equal}(\pi_2,\tp_2)^{-1}(g,h,\xi)\in E_2$.
By Equation~(\ref{eqn:phi-n}) above, $\phi_2(x)
{\equal}(g,h,g(hx),hx,x)$. It follows that $\phi_1(\teps_1(x))
{\equal}\teps_1(\phi_2(x)){\equal}(gh,g(hx),x)$. Using again
(\ref{eqn:action}), we get $(gh)x{\equal}g(hx)$.

To show that $p(\xi)\xi{\equal}\xi$ for all $\xi\in E_0$, let
$x{\equal}\teta_0(\xi)\in E_1$, then
$$\phi_1(x){\equal}(\pi_1(x),\tp_0(x),\tp_1(x)){\equal}(p(\xi),\xi,\xi),$$
therefore, by (\ref{eqn:action}), we get $p(\xi)\xi{\equal}\xi$.

We have shown that $G$ acts on $E_0$. To show that the simplicial
spaces $(G\ltimes E_0)_\com$ and $E_\com$ are isomorphic,
let $\phi'_n\colon G_n\times_{\tp_n,p}E_0{\equal}(G\ltimes E_0)_n
\to G_n\times E_0^{n+1}$ be defined like $\phi$, then $\phi'_n$
is also injective for all $n$ and $\phi'_\com$ is a morphism
of simplicial spaces. By Equation~(\ref{eqn:phi-n}), we have
a commutative diagram
$$\xymatrix{
E_n \ar[r]^{\hspace{-5ex}(\pi_n,\tp_n)}\ar[dr]_{\phi_n}
& \quad G_n\times_{\tp_n,p}E_0\ar[d]^{\phi'_n}\\
&G_n\times E_0^{n+1}
}$$
where $\phi_n$ and $\phi'_n$ are injective and
$(\pi_n,\tp_n)$ is a homeomorphism. Therefore,
$(\pi_n,\tp_n)_{n\in \N}$ is an isomorphism of simplicial spaces.
\pfend

\section{\v Cech cohomology}\label{sec:cech}
\subsection{Covers of simplicial spaces}
\begin{defi}
An \emph{open cover} of a semi-simplicial space $M_\com$ is a sequence
$\uU_\com{\equal}(\uU_n)_{n\in \N}$ such that $\uU_n{\equal}(U^n_i)_{i\in I_n}$
is an open cover of the space $M_n$.

The cover is said to be \emph{semi-simplicial} if $I_\com
{\equal}(I_n)_{n\in \N}$ is a semi-simplicial set such that for all
$f\in\hom_{\Delta}(k,n)$ and for all $i\in I_n$ one has
$\tf(U^n_i){\equal}U^k_{\tf(i)}$. In the same way, one defines the
notions of simplicial cover and of $N$-simplicial cover.
\end{defi}

The reason why we need to introduce this terminology is that,
even when $M_\com$ is a simplicial space, there may not exist
sufficiently fine simplicial covers. However, given a cover $\uU_\com$,
we can form the semi-simplicial cover $\sigma\uU_\com$
defined as follows.

Let $\pP_n{\equal}\cup_{k{\equal}0}^n \pP^k_n$, where
$\pP_n^k{\equal}\hom_{\Delta'}(k,n)$.
Note that $\pP_n$ can be identified with the set of nonempty
subsets of $[n]$.

Let $\Lambda_n$ (or $\Lambda_n(I)$ if there is a risk of confusion)
be the set of maps
\begin{equation}\label{eqn:Lambda}
\lambda\colon \pP \to \cup_k I_k\,\mbox{ such that }
\lambda(\pP_n^k)\subset I_k.
\end{equation}
For all $\lambda\in \Lambda_n$,
we let
$$U^n_\lambda{\equal}\bigcap_{k\le n}\bigcap_{f\in \pP_n^k}
\tf^{-1}(U^k_{\lambda(f)}).$$
It is clear that
$(U_\lambda)_{\lambda\in \Lambda_n}$ is an open cover of $M_n$.
\par
The semi-simplicial structure on $\Lambda_\com$ is defined
in an obvious way: for all $g\in\hom_{\Delta'}(n,n')$,
$\tg\colon \Lambda_{n'}\to\Lambda_n$ is the map
$$(\tg\lambda')(f){\equal}\lambda'(g\circ f).$$

In the same way, for all integers $n\le N$, let
\begin{equation}\label{eqn:sigmaN}
(\sigma_N\uU)_n{\equal}(U^n_\lambda)_{\lambda\in \Lambda^N_n},
\end{equation}
where $\Lambda_n^N$ is the set of all maps
$\displaystyle\lambda\colon \bigcup_{k\le n}\hom_\Delta(k,n)\to
\bigcup_{k\le n}I_k$ which satisfy
$\lambda(\hom_\Delta(k,n))\subset I_k$, and
$$U^n_\lambda{\equal}\bigcap_{k\le n}\bigcap_{f\in\hom_\Delta(k,n)}
\tf^{-1}(U^k_\lambda(f)).$$

The $N$-simplicial structure on $\Lambda^N_\com$ is defined
as follows: for all integers $n,n'\le N$ and all
$g\in \hom_\Delta(n,n')$, $\tg\colon \Lambda^N_{n'}
\to\Lambda^N_n$ is the map
$(\tg\lambda')(f){\equal}\lambda'(g\circ f)$.

Then $\sigma_N\uU_\com{\equal}(\sigma_N\uU_n)_{n\le N}$ is a
$N$-simplicial cover of the $N$-skeleton of $M_\com$.

\begin{convention}\label{conv:sigmaN}
We will also (abusively) denote by $\sigma_N\uUb$ the open cover
which coincides with $\sigma_N\uUb$ for $n\le N$ and with $\uUb$
for $n\ge N+1$.
\end{convention}

\begin{example}\label{ex:constant}
Let $M_\com{\equal}(M)_{n\in \N}$ be the constant simplicial space
associated to a topological space $M$, and
suppose $\uU_0{\equal}(U^0_i)_{i\in I_0}$ is an open cover of $M$.
Define $I_n{\equal}I_0^{n+1}$, then $I_\com{\equal}(I_n)_{n\in\N}$ is
endowed with a simplicial structure by $\tf(i_0,\ldots,i_n)
{\equal}(i_{f(0)},\ldots,i_{f(k)})$ for all $f\in\hom_\Delta(k,n)$.
Let $U^n_{(i_0,\ldots,i_n)}{\equal}U^0_{i_0}\cap\cdots\cap
U^0_{i_n}$, and let $\uU_n{\equal}(U_i^n)_{i\in I_n}$, then
$\uUb$ is a simplicial cover of $M_\com$.
\end{example}

The ``set'' of covers of a simplicial space $M_\com$ is endowed
with a partial preorder. Suppose $\uU_\com$ and $\vV_\com$
are open covers of $M_\com$, with $\uU_n
{\equal}(U^n_i)_{i\in I_n}$ and $\vV_n{\equal}(V^n_j)_{j\in J_n}$. We say
that $\vV$ is \emph{finer} than $\uU$ if for all $n$ there exists
$\theta_n\colon J_n\to I_n$ such that $\theta_n(V^n_j)
\subset U^n_{\theta(j)}$ for all $j$. The map
$\theta_\com{\equal}(\theta_n)_{n\in\N}$ is required to be
semi-simplicial (resp. $N$-simplicial) if $\uU$ and $\vV$
are semi-simplicial (resp. $N$-simplicial).

\subsection{\v{C}ech cohomology}
Let $\uU_\com$ be a semi-simplicial open cover of $M_\com$
and let $\aA^\com$ be a semi-simplicial abelian sheaf.
Define a complex
$$C^n_{ss}(\uU_\com;\aA^\com)
{\equal}\prod_{i\in I_n}\aA^n(U^n_i),$$
i.e. $C^n_{ss}(\uU_\com;
\aA^\com)$ is the space of global sections of the pull-back
of $\aA^n$ on $\coprod_{i\in I_n} U^n_i$.

Define the differential $d\colon C_{ss}^n(\uU_\com;\aA^\com)
\to C_{ss}^{n+1}(\uU_\com;\aA^\com)$ by
$$(dc)_i{\equal}\sum_{k{\equal}0}^{n+1} (-1)^k \teps^*_k c_{\teps_k(i)},$$
where $\teps^*_k c_{\teps_k(i)}$ is the ``restriction''
of $c_{\teps_k(i)}\in \aA^n(U^n_{\teps_k(i)})$ to a section in
$\aA^{n+1}(U^{n+1}_i)$.

It is immediate to check that $d^2{\equal}0$, hence we may define
the cohomology groups $H^*_{ss}(\uU_\com;\aA^\com)$.

\begin{example}
In Example~\ref{ex:constant}, suppose $\aA$ is an abelian sheaf
on $M$ and that $\aA^n{\equal}\aA$ for all $n$. Then $H^*_{ss}(\uU_\com;
\aA^\com)$ is identical to the usual cohomology group
$H^*(\uU_0;\aA)$.
\end{example}

Let $\uU_\com$ be any open cover of $M_\com$. We denote
\begin{eqnarray}\label{eqn:definition-Cn(U)}
C^n(\uU_\com;\aA^\com) &{\equal}& C^n_{ss}(\sigma\uU_\com;
\aA^\com)\\
\label{eqn:definition-Hn(U)}
H^n(\uU_\com;\aA^\com) &{\equal}& H^n_{ss}(\sigma\uU_\com;
\aA^\com).
\end{eqnarray}

Now, we want to define \v{C}ech cohomology. The idea is to define
$\vH^n(M_\com;\aA^\com)$ as the inductive limit
over $\uU_\com$ of the groups $H^n(\uU_\com;\aA^\com)$.
The problem is that if $\theta_\com\colon J_\com
\to I_\com$ is a refinement, then $\theta_\com$
indeed defines a restriction map

\begin{eqnarray}
\theta^*\colon C^n(\uU_\com;\aA^\com)
&\to& C^n(\vV_\com;\aA^\com)\\\label{eqn:definition-theta*}
(\theta^*\varphi)_j&{\equal}&\mbox{restriction to }V^n_j\mbox{ of }
\varphi_{\theta_n(j)}
\end{eqnarray}

which commutes with the differentials,
and thus $\theta^*$ defines a map
$$\theta^*\colon H^n(\uU_\com;\aA^\com)
\to H^n(\vV_\com;\aA^\com).$$
However, that map may depend on
the choice of $\theta$. On the other hand we have the

\begin{lem}\label{lem:dH+Hd}
Let $N\in\N$. Suppose that $\uU_\com$ and $\vV_\com$
are open covers of $M_\com$ such that $\vVb$
admits a $N$-simplicial structure. Suppose that $\vV_\com$
is finer than $\uU_\com$ and that $\theta_0$,
$\theta_1\colon \uU_\com\to \vV_\com$ are two refinements.
Then for all $n\le N$
there exists $H\colon C^n(\uU_\com;\aA^\com)
\to C^{n-1}(\vV_\com;\aA^\com)$ (with the convention
$C^{-1}{\equal}\{0\}$) such that $\theta_1^*-\theta_0^*{\equal}dH+Hd$.
\end{lem}

(In the lemma above, we say that $\vVb$ has a $N$-simplicial structure
if the $N$-skeleton $(\vV_n)_{n\in\N}$ has a $N$-simplicial structure.)

\pf
Define for all $\varphi\in C^n(\uU_\com;\aA^\com)$ and
for all $\lambda\in \Lambda_{n-1}(J)$ (recall notation (\ref{eqn:Lambda})):
\begin{equation}
\label{eqn:Hphi}
(H\varphi)_\lambda {\equal} \sum_{k{\equal}0}^{n-1} (-1)^k
\teta_k^*\varphi_{\alpha_k(\lambda)},
\end{equation}
where as usual $\eta_k\colon [n]\to [n-1]$ is the $k$-th
degeneracy map, and $\alpha_k$ is defined as follows: for
all $f\in \hom_{\Delta'}(r,n)$, let
$$\alpha_k(\lambda)(f){\equal}
\left\{
\begin{array}{l}
\theta_0(\lambda(\eta_k\circ f))\;\mbox{ if }
  \{k,k+1\}\not\subset f([r])\mbox{ and }f(0)\le k\\
\theta_1(\lambda(\eta_k\circ f))\;\mbox{ if }
  \{k,k+1\}\not\subset f([r])\mbox{ and }f(0)\ge k+1\\
\theta_0(\teta_{k'}(\lambda(f'))\;\mbox{ if }
  \{k,k+1\}\subset f([r])
\end{array}
\right.$$
where in the third line,
$k'$ is the integer such that $f(k'){\equal}k$ and $f'$ is the unique
morphism in $\hom_{\Delta'}(r-1,n-1)$ such that the following
diagram commutes:
\begin{equation}\label{eqn:diag-rfn}
\xymatrix{
[r]\ar[r]^f\ar[d]^{\eta_{k'}} & [n]\ar[d]^{\eta_k}\\
[r-1]\ar[r]^{f'} & [n-1]
}
\end{equation}
i.e. $f'(i){\equal}f(i)$ for $i\le k'$ and $f'(i){\equal}f(i+1)-1$ for
$i\ge k'+1$.

Let us first check that formula (\ref{eqn:Hphi}) makes sense, i.e. that
$V^{n-1}_\lambda \subset \teta_k^{-1}(U^r_{\alpha_k(\lambda)})$
for all $k$.

Since $\displaystyle U^n_{\alpha_k(\lambda)}
{\equal}\bigcap_{r{\equal}0}^n \bigcap_{f\in\hom_{\Delta'}(r,n)}
\tf^{-1}(U^n_{\alpha_k(\lambda)(f)})$, we need to show that
\begin{equation}\label{eqn:inclusion}
V^{n-1}_\lambda \subset \teta_k^{-1}\tf^{-1}(U^r_{\alpha_k(\lambda)(f)}).
\end{equation}
If $\{ k,k+1 \}\not\subset f([r])$ and $f(0)\le k$ then
\begin{eqnarray*}
V^{n-1}_\lambda &\subset& (\widetilde{\eta_k\circ f})^{-1}
(V^r_{\lambda(\eta_k\circ f)})\quad\mbox{by definition of $V^{n-1}_\lambda$}\\
&{\equal}&\teta_k^{-1}\tf^{-1}(V^r_{\lambda(\eta_k\circ f)})\\
&\subset& \teta_k^{-1}\tf^{-1}(U^r_{\theta_0(\lambda(\eta_k\circ f))})
\quad\mbox{since $\theta_0\colon J_\com\to I_\com$
is a refinement}\\
&{\equal}&\teta_k^{-1}\tf^{-1}(U^r_{\alpha_k(\lambda)(f)}).
\end{eqnarray*}
\par\smallskip
If $\{k,k+1\}\not\subset f([r])$ and $f(0)\ge k+1$ the proof of
(\ref{eqn:inclusion}) is the same,
except that $\theta_0$ is replaced by $\theta_1$.
\par\smallskip
If $\{k,k+1\}\subset f([r])$ then
\begin{eqnarray*}
V^{n-1}_\lambda &\subset& (\widetilde{f'})^{-1}(V^{r-1}_{\lambda(f')})\\
&\subset& \widetilde{f'}^{-1}\teta_{k'}^{-1}(V^r_{\teta_k(\lambda(f'))})
\quad\mbox{(recall $\vVb$ is $N$-simplicial)}\\
&\subset& \teta^{-1}_k\tf^{-1}(U^r_{\theta_0(\teta_k(\lambda(f')))})
\quad\mbox{by (\ref{eqn:diag-rfn})},
\end{eqnarray*}
thus (\ref{eqn:inclusion}) is proved.
\par\bigskip
Let us show that
\begin{equation}\label{eqn:dH+Hd}
dH+Hd {\equal} \theta_1^* - \theta_0^*.
\end{equation}
We have
$$(Hd\varphi)_\lambda{\equal}\sum_{\ell{\equal}0}^n\sum_{k{\equal}0}^{n+1} A_{k,\ell}
\;\mbox{ and }\; (dH\varphi)_\lambda{\equal}\sum_{k{\equal}0}^n\sum_{\ell{\equal}0}^{n-1}
B_{k,\ell},$$
where $A_{k,\ell}{\equal}(-1)^{k+\ell}\teta^*_\ell\teps^*_k
\varphi_{\teps_k(\alpha_\ell(\lambda))}$ and
$B_{k,\ell}{\equal}(-1)^{k+\ell} \teps_k^*\teta^*_\ell
\varphi_{\alpha_\ell(\teps_k(\lambda))}$.

We have $A_{\ell,\ell}{\equal}\teta_\ell^*\teps_\ell^*
\varphi_{\teps_\ell(\alpha_\ell(\lambda))} {\equal}
\varphi_{\teps_\ell(\alpha_\ell(\lambda))}$, and for all $f\in \pP_n$,
$$\teps_\ell(\alpha_\ell(\lambda))(f)
{\equal}\alpha_\ell(\lambda)(\varepsilon_\ell\circ f)
{\equal}\theta_j(\lambda(\eta_\ell\circ\varepsilon_\ell(f)))
{\equal}\theta_j(\lambda(f)),$$
where $j{\equal}0\iff \varepsilon_\ell\circ f(0)\le\ell
\iff f(0) \le \ell-1$, and $j{\equal}1$ otherwise.

Let $\lambda^{(p)}(f) {\equal} \left\{
\begin{array}{ll}
\theta_0(\lambda(f))&\mbox{if } f(0)\le p\\
\theta_1(\lambda(f))&\mbox{if } f(0)\ge p+1,
\end{array}
\right.$
then
\begin{equation}\label{eqn:All}
A_{\ell,\ell}{\equal}\varphi_{\lambda^{\ell-1}}.
\end{equation}
Similarly, we have $A_{\ell+1,\ell}{\equal}-\teta_\ell^*\teps_{\ell+1}^*
\varphi_{\teps_{\ell+1}(\alpha_\ell(\lambda))} {\equal}
-\varphi_{\teps_{\ell+1}(\alpha_\ell(\lambda))}$, and for all
$f\in \pP_n$,
$$\teps_{\ell+1}(\alpha_\ell(\lambda))(f)
{\equal}\alpha_\ell(\lambda)(\varepsilon_{\ell+1}\circ f)
{\equal}\theta_j(\lambda(\eta_\ell\circ\varepsilon_{\ell+1}\circ f))
{\equal}\theta_j(\lambda(f)),$$
where $j{\equal}0\iff \varepsilon_{\ell+1}\circ f(0)\le\ell \iff
f(0)\le \ell$, and $j{\equal}1$ otherwise. We thus get
\begin{equation}\label{eqn:Al+1l}
A_{\ell+1,\ell}{\equal}-\varphi_{\lambda^{\ell}}.
\end{equation}
From (\ref{eqn:All}) and (\ref{eqn:Al+1l}) we obtain
$\sum_{\ell\le k\le \ell+1} A_{k,\ell} {\equal}\varphi_{\lambda^{(-1)}}
-\varphi_{\lambda^{(n)}}
{\equal}\theta^*_1\varphi-\theta_0^*\varphi$.
\par\smallskip
Let us examine the other terms. We have $\sum_{k\le\ell-1} A_{k,\ell}
{\equal}\sum_{0\le k\le \ell\le n-1}A_{k,\ell+1}$ and
$\sum_{k\ge\ell+2} A_{k,\ell}
{\equal}\sum_{0\le\ell \inferieur k\le n} A_{k+1,\ell}$.
To complete the proof of (\ref{eqn:dH+Hd}) it suffices to show that
$A_{k,\ell+1}+B_{k,\ell}{\equal}0$ for $k\le\ell$, and that
$A_{k+1,\ell}+B_{k,\ell}{\equal}0$ for $k\ge \ell+1$. Noting that
$\eta_{\ell+1}\varepsilon_k {\equal} \varepsilon_k\eta_\ell$ for $k\le\ell$
and that $\eta_\ell\varepsilon_{k+1}{\equal}\varepsilon_k\eta_\ell$ for
$k\ge \ell+1$, it suffices to show that
\begin{itemize}
\item[(a)] $\teps_k(\alpha_{\ell+1}(\lambda))
{\equal}\alpha_\ell(\teps_k(\lambda))$ for $k\le\ell$,
\item[(b)] $\teps_{k+1}(\alpha_\ell(\lambda)){\equal}\alpha_\ell
(\teps_k(\lambda))$ for $k\ge \ell+1$.
\end{itemize}

Let us show (a). Suppose that $f\in\hom_{\Delta'}(r,n)$
and let us first treat the case
$\{\ell,\ell+1\}\not\subset f([r])$. Then, letting $j{\equal}0$
for $\varepsilon_k\circ f(0)\le\ell+1$ ($\iff f(0) \le\ell$),
and $j{\equal}1$ otherwise, we have
\begin{eqnarray*}
\lefteqn{\teps_k(\alpha_{\ell+1}(\lambda))(f)
{\equal}\alpha_{\ell+1}(\lambda)(\varepsilon_k\circ f)
{\equal}\theta_j(\lambda(\eta_{\ell+1}\circ\varepsilon_k\circ f))}\\
&{\equal}&\theta_j(\lambda(\varepsilon_k\circ\eta_\ell\circ f))
{\equal}\alpha_\ell(\teps_k(\lambda))(f).
\end{eqnarray*}

Let us treat the case $\{\ell,\ell+1\}\subset f([r])$.
Let $\ell'$ such that $f(\ell'){\equal}\ell$, and let
$f'\colon [r-1]\to [n-1]$ be the increasing map such that the diagram
\begin{equation}\label{eqn:diag-rfn2}
\xymatrix{
{[r]}\ar[r]^{f}\ar[d]^{\eta_{\ell'}} &[n]\ar[d]^{\eta_\ell}\\
{[r-1]}\ar[r]^{f'} & [n-1]
}
\end{equation}
commutes. Since $\varepsilon_k\circ\eta_\ell {\equal} \eta_{\ell+1}\circ
\varepsilon_k\colon [n]\to [n]$, the diagram
\begin{equation}\label{eqn:diag-rfn+1}
\xymatrix{
{[r]}\ar[r]^{\varepsilon_k\circ f}\ar[d]^{\eta_{\ell'}}
&[n+1]\ar[d]^{\eta_{\ell+1}}\\
{[r-1]}\ar[r]^{\varepsilon_k\circ f'} & [n]
}
\end{equation}
commutes. We thus see that
\begin{eqnarray*}
\teps_k(\alpha_{\ell+1}(\lambda))(f)&{\equal}&\alpha_{\ell+1}(\lambda)
(\varepsilon_k\circ f)\\
&{\equal}&\theta_0(\teta_{\ell'}(\lambda(\varepsilon_k
\circ f')))\quad\mbox{by }(\ref{eqn:diag-rfn+1})\\
&{\equal}&\theta_0(\teta_{\ell'}((\teps_k\lambda)(f')))\\
&{\equal}&\alpha_\ell(\teps_k\lambda)(f)\quad\mbox{by }
(\ref{eqn:diag-rfn2}).
\end{eqnarray*}
This completes the proof of (a).
\par\medskip
Let us show (b): the method is similar. If $\{\ell,\ell+1\}
\not\subset f([r])$ then
\begin{eqnarray*}
\lefteqn{\teps_{k+1}(\alpha_\ell(\lambda))(f)
{\equal}\alpha_\ell(\lambda)(\varepsilon_{k+1}\circ f)
{\equal}\theta_j(\lambda(\eta_\ell\circ\varepsilon_{k+1}\circ f))}\\
&{\equal}&\theta_j(\lambda(\varepsilon_k\circ\eta_\ell\circ f))
{\equal}\theta_j((\teps_k\lambda)(\eta_\ell\circ f))
{\equal}\alpha_\ell(\teps_k(\lambda))(f).
\end{eqnarray*}
If $\{\ell,\ell+1\}\subset f([r])$, let $\ell'$ such that
$f(\ell'){\equal}\ell$ and let $f'$ such that
\begin{equation}\label{eqn:diag-rfn3}
\xymatrix{
{[r]}\ar[r]^{f}\ar[d]^{\eta_{\ell'}}
&[n]\ar[d]^{\eta_\ell}\\
{[r-1]}\ar[r]^{f'}&[n-1]
}
\end{equation}
commutes. Then since $\varepsilon_k\circ\eta_\ell{\equal}\eta_\ell
\circ\varepsilon_{k+1}\colon [n]\to [n]$, the diagram
\begin{equation}\label{eqn:diag-rfn+12}
\xymatrix{
{[r]}\ar[r]^{\varepsilon_{k+1}\circ f}\ar[d]^{\eta_{\ell'}}
&[n+1]\ar[d]^{\eta_\ell}\\
{[r-1]}\ar[r]^{\varepsilon_k\circ f'}&[n]
}
\end{equation}
commutes, therefore
\begin{eqnarray*}
\teps_{k+1}(\alpha_\ell(\lambda))(f)
&{\equal}&\alpha_\ell(\lambda)(\varepsilon_{k+1}\circ f)\\
&{\equal}&\theta_0(\teta_{\ell'}(\lambda(\varepsilon_k
\circ f')))\quad\mbox{by }(\ref{eqn:diag-rfn+12})\\
&{\equal}&\theta_0(\teta_{\ell'}((\teps_k\lambda)(f')))\\
&{\equal}&(\alpha_\ell(\teps_k\lambda))(f)\quad
\mbox{by }(\ref{eqn:diag-rfn3}).
\end{eqnarray*}
This completes the proof of (b) and hence of (\ref{eqn:dH+Hd}).
\pfend

Let us now define
\begin{equation}\label{eqn:definition-cech}
\vH^n(M_\com;\aA^\com){\equal}
\lim_\to H^n(\uU_\com;A^\com)
\end{equation}
where $\uU_\com$ runs over open covers of $M_\com$
whose $N$-skeleton admits a $N$-simplicial structure for
some $N\ge n+1$. (Recall that $H^n(\uU_\com;\aA^\com)$
was defined by equation (\ref{eqn:definition-Hn(U)}).)

To avoid set-theoretic difficulties (since the collection of open covers
is not a set), we can restrict ourselves to open covers indexed
by sets of cardinality $\le\sum_n\# M_n$.

By Lemma~\ref{lem:dH+Hd} above, if $\vV_\com$ is finer than $\uU_\com$
and if $\vV_\com$ has a $N$-simplicial structure, then there is a
canonical map $H^*(\uUb;\aAb)\to H^*(\vVb;\aAb)$ defined by
Eqn.~(\ref{eqn:definition-theta*}) ($\theta$ is \emph{not} required
to respect the $N$-simplicial structures).
Since every open cover of $\Mb$ admits a $N$-simplicial refinement
(see Convention~\ref{conv:sigmaN}),
the inductive limit is well-defined and is an abelian group.

Moreover, for every open cover $\uUb$ of $\Mb$, $N$-simplicial or not,
there is a canonical map $H^n(\uUb;\aAb)\to \vH^n(\Mb;\aAb)$
obtained by mapping $H^n(\uUb;\aAb)$ to $H^n(\vVb;\aAb)$
using (\ref{eqn:definition-theta*}), where $\vVb$ is any refinement
admitting a $N$-simplicial structure for some $N\ge n+1$.

It is clear that any element $[\varphi]$ of $H^n(\uUb;\aAb)$ maps
to 0 in $\vH^n(\Mb;\aAb)$ if and only if there exists a refinement
($N$-simplicial or not) such that $[\varphi]$ maps to 0
in $H^n(\vVb;\aAb)$. Thus, in some sense, we can say that
$\vH^n(\Mb,\aAb)$ is the inductive limit of $H^n(\uUb,\aAb)$
where $\uUb$ runs over \emph{all} open covers of $\Mb$.

\begin{example}
Consider the case of a discrete group $G$, and suppose that
$\aAb$ is the sheaf associated to a $G$-module $A$
(Definition~\ref{def:G-module}). Then, from (\ref{eqn:simplicial-GxA})
and below, we see that
\begin{eqnarray*}
\lefteqn{(dc)_\lambda(g_1,\ldots,g_{n+1})
{\equal}g_1c_{\teps_0\lambda}(g_2,\ldots,g_{n+1})+}\\
&&+\sum_{k{\equal}1}^n (-1)^kc_{\teps_k\lambda}(g_1,\ldots,g_kg_{k+1},
\ldots,g_{n+1})+\\
&&\qquad + (-1)^{n+1}c_{\teps_{n+1}\lambda}
(g_1,\ldots,g_n).
\end{eqnarray*}

(Compare with (\ref{eqn:moore-differential}).)
Considering the maximal open cover $(U^n_x)_{x\in G^n}$ where
$U^n_x{\equal}\{x\}$, one easily sees that \v{C}ech cohomology coincides
with usual group cohomology.
\end{example}

\begin{rem}\label{rem:cech-complex}
As in \cite{God73}, $\vH^n(\Mb;\aAb)$ can be seen as the $n$-th
cohomology group of a canonical \v{C}ech complex. Indeed, let
$\mathcal{R}(\Mb)$ be the set of covers of the form $\uU_n
{\equal}(U^n_x)_{x\in M_n}$. If $\uUb$ and $\vVb$ are in $\mathcal{R}(\Mb)$,
let us say that $\vVb$ is finer than $\uUb$ if for all $n$
and all $x\in M_n$ we have $V^n_x\subset U^n_x$.
Given $\uUb$ in $\mathcal{R}(\Mb)$, denote by $\sigma_N\uUb$
the associated $N$-simplicial cover (see Convention~\ref{conv:sigmaN})
and let
$$\check{C}^n_N(\Mb;\aAb){\equal}\lim_{\uUb} C^n(\sigma_N\uUb;\aAb)
:{\equal}\lim_{\uUb} C^n_{ss}(\sigma\sigma_N\uUb,\aAb),$$
where $\uUb$ runs over open covers in $\mathcal{R}(\Mb)$.

Then $\vH^n(\Mb;\aAb)$ is the cohomology of $\check{C}^n_N
(\Mb;\aAb)$ whenever $N\ge n+1$.
\end{rem}

\subsection{Compatibility with usual \v{C}ech cohomology for
spaces}

Let $M$ be a space and $\aA$ an abelian sheaf on $M$. Denote by
$\Mb$ be the constant simplicial space associated to $M$ and by
$\aAb$ the sheaf on $\Mb$ corresponding to $\aA$.

We want to show the

\begin{prop}
With the above assumptions, the usual \v{C}ech cohomology groups
$\vH^n(M;\aA)$ are isomorphic to the \v{C}ech cohomology groups
$\vH^n(\Mb;\aAb)$.
\end{prop}

\pf
To determine $\vH^n(\Mb;\aAb)$ we can restrict ourselves to covers
$\uUb$ of the form $\uU_n{\equal}(U^n_i)_{i\in I_n}$ where $I_n{\equal}I_0^{n+1}$
and $U^n_{i_0,\ldots,i_n}{\equal}U^0_{i_0}\cap\cdots\cap U^0_{i_n}$.
Let us show that $H^*(\uUb;\aAb)\cong H^*(\uU_0;\aA)$. It is not
obvious that these two groups are isomorphic, since
$C^*(\uUb;\aAb){\equal}C_{ss}^*(\sigma\uUb;\aAb)$ while
$C^*(\uU_0;\aA){\equal}C^*_{ss}(\uUb;\aAb)$. However, we show that these
two complexes are homotopically equivalent:

First, there is an obvious map
$$q\colon C^*(\uUb;\aAb)\to C^*(\uU_0;\aA)$$
defined by $(q\varphi)_{i_0,\ldots,i_n}{\equal}\varphi_{\lambda^{(i)}}$,
where
$$\lambda^{(i)}(f){\equal}(i_{f(0)},\ldots,i_{f(r)})\quad\mbox{for all }
f\in\hom_{\Delta'}(r,n).$$
In the other direction, define $\iota\colon C^*(\uU_0,\aA)
\to C^*(\uUb;\aAb)$ by
$$(\iota c)_\lambda {\equal} c_{\lambda_0,\ldots,\lambda_n},$$
where $\lambda_k$ denotes $\lambda(p_k)$ and $p_k\colon [0]\to [n]$
denotes the map $p_k(0){\equal}k$.

We have $q\circ \iota{\equal}\id$. Indeed, $((q\circ\iota)(c))_{i_0,\ldots,i_n}
{\equal}(\iota c)_{\lambda^{(i)}} {\equal} c_{i_0,\ldots,i_n}$.

Conversely, we don't have $\iota\circ q{\equal}\id$ since
$((\iota\circ q)(\varphi))_\lambda{\equal}\varphi_{\lambda'}$, where
$$\lambda'(f){\equal}(\lambda_{f(0)},\ldots,\lambda_{f(r)})
\mbox{ for all }f\in\hom_{\Delta'}(r,n).$$
However, $\iota\circ q$ and $\id$
are homotopic. Indeed, define $H\colon C^n(\uUb;\aAb)\to
C^{n-1}(\uUb;\aAb)$ by
\begin{equation}\label{eqn:H}
(H\varphi)_\lambda {\equal} \sum_{k{\equal}0}^{n-1} (-1)^k
\teta_k^*\varphi_{\alpha_k(\lambda)}
\end{equation}
where 
$$\alpha_k(\lambda)(f){\equal}
\left\{
\begin{array}{l}
\lambda(\eta_k\circ f)\;\mbox{ if }
  \{k,k+1\}\not\subset f([r])\mbox{ and }f(0)\le k\\
\lambda(\eta_k\circ f)\;\mbox{ if }
  \{k,k+1\}\not\subset f([r])\mbox{ and }f(0)\ge k+1\\
\teta_{k'}(\lambda(f'))\;\mbox{ if }
  \{k,k+1\}\subset f([r])
\end{array}
\right.$$
($f(k'){\equal}k$ and $f'$ is defined as in the proof of
Lemma~\ref{lem:dH+Hd}; also, recall that $\teta_{k'}(i_0,\ldots,i_{r-1})
{\equal}(i_0,\ldots,i_{k'},i_{k'},\ldots,i_{r-1})$.)

Then the same proof as in Lemma~\ref{lem:dH+Hd}
shows that $dH+Hd{\equal}\iota\circ q - \id$.
We leave out details; anyway we will show later that \v{C}ech
cohomology coincides with sheaf cohomology for paracompact simplicial
spaces, so (at least in the paracompact case) this will provide a second proof
that sheaf cohomology for spaces coincides with sheaf cohomology for
the associated constant simplicial space.
\pfend

Let us introduce some notation:

\begin{notation}\label{not:cocycle}
For any simplicial space $\Mb$ and any open cover $\uUb$, let us write
elements $\lambda\in\Lambda_n$ (see (\ref{eqn:Lambda})) as
($2^{n+1}-1$)-tuples $(\lambda_S)_{\emptyset\ne S\subset [n]}$,
where subsets $S$ are ordered first by cardinality, then by lexicographic
order. For instance, the triple $(\lambda_0,\lambda_1,\lambda_{01})$
represents the element $\lambda\in \Lambda_1$ such that
$\lambda(\{0\}){\equal}\lambda_0$, $\lambda(\{1\}){\equal}\lambda_1$,
$\lambda(\{0,1\}){\equal}\lambda_{01}$. A cochain in $C^1(\uUb;\aAb)$
is thus a family $(\varphi_{\lambda_0\lambda_1\lambda_{01}})$.
\end{notation}

Then, we can write (\ref{eqn:H}) more explicitly. For instance,
the formulas for $n{\equal}1$ and $n{\equal}2$ are respectively

\begin{eqnarray*}
(H\varphi)_{\lambda_0}&{\equal}&\varphi_{\lambda_0\lambda_0\lambda'_{00}}\\
(H\varphi)_{\lambda_0\lambda_1\lambda_{01}}&{\equal}&
\varphi_{\lambda_0\lambda_0\lambda_1\lambda'_{00}\lambda_{01}
\lambda'_{01}\lambda'_{001}}
-\varphi_{\lambda_0\lambda_1\lambda_1\lambda_{01}\lambda_{01}
\lambda'_{11}\lambda'_{011}},
\end{eqnarray*}
where $\lambda'_{i_0\cdots i_r}$ denotes the $r+1$-tuple
$(\lambda_{i_0},\ldots,\lambda_{i_r})\in I_0^{r+1}{\equal}I_r$.

\subsection{Long exact sequences in \v{C}ech cohomology}
In this section, most proofs are almost identical to \cite{God73},
thus we will only sketch them.

\begin{prop}\label{prop:presheaf-exact}
If $0\to {\aA'}^\com \to \aAb \to {\aA''}^\com\to 0$ is an exact
sequence of abelian \underline{presheaves} then the functor $\aA\mapsto
\vC^*_N(\Mb;\aAb)$ (see Remark~\ref{rem:cech-complex})
maps the above exact sequence to an exact sequence of complexes.
\end{prop}

\pf
$0\to C^*_{ss}(\sigma\sigma_N\uUb;{\aA'}^\com)
 \to C^*_{ss}(\sigma\sigma_N\uUb;\aAb)
 \to C^*_{ss}(\sigma\sigma_N\uUb;{\aA''}^\com)\to 0$
is exact for every open cover $\uUb$.
\pfend

Let us say that a simplicial space $\Mb$ is \emph{paracompact} if each
$M_n$ is paracompact.

\begin{prop}\label{prop:zero-sheaf}\cite[theorem~5.10.2]{God73}
If $\aAb$ is an abelian presheaf on a paracompact simplicial space $\Mb$
such that $\aAb$ induces the zero sheaf then $\vH^n(\Mb;\aAb){\equal}0$
for all $n\ge 0$.
\end{prop}

\pf
Using paracompactness,
every cohomology class is represented by a cocycle in
$C^n_N(\uUb;\aAb)$ with $\uU_n$ locally finite $\forall n$.
Then, using the fact that each $\aA^n$ induces the zero sheaf,
every {\em{cochain}} of that cover becomes zero once restricted
to a suitable finer cover.
\pfend

\begin{coro}\label{cor:presheaf-sheaf}
If an abelian presheaf $\aAb$ over a paracompact simplicial space
$\Mb$ induces the sheaf $\tilde{\aA}^\com$ then
$\vH^n(\Mb,\aAb)\stackrel{\simeq}{\to}
\vH^n(\Mb;\tilde{\aA}^\com)$.
\end{coro}

\pf
There are exact sequences of presheaves
\begin{eqnarray*}
0\to \nN^\com \to \aAb \to \jJ^\com\to 0\\
0\to \jJ^\com \to \tilde{\aA}^\com \to \qQ^\com\to 0,
\end{eqnarray*}
where $\nN^\com$ and $\qQ^\com$ induce the zero sheaf:
$\nN^n(U)$ is the set of sections in $\aA^n(U)$ whose germ at every point
is zero, and $\tilde{\aA}^n(U){\equal}\{(\sigma_i)_{i\in I}\}/\sim$,
where $\sigma_i\in \jJ^n(U_i)$ for some open cover $(U_i)_{i\in I}$ of
$U$ and the equivalence relation $\sim$ is defined by
$(\sigma_i)_{i\in I}\sim (\sigma'_j)_{j\in J}$ iff
$\forall i,j$, ${\sigma_i}_{|U_i\cap U'_j}{\equal}
{\sigma'_j}_{|U_i\cap U'_j}$.
The conclusion follows from Propositions~\ref{prop:presheaf-exact}
and~\ref{prop:zero-sheaf} above.
\pfend

\begin{coro}\label{cor:long-exact-cech}
If $0\to {\aA'}^\com \to \aAb \to {\aA''}^\com$ is an exact
sequence of sheaves over a paracompact simplicial space $\Mb$,
then there is a natural long exact sequence
$$0\to \vH^0(\Mb;{\aA'}^\com) \to
\vH^0(\Mb; \aAb) \to 
\vH^0(\Mb;{\aA''}^\com)
\stackrel{\partial}{\to}\vH^1(\Mb;{\aA'}^\com)\to\cdots$$
\end{coro}

\pf
Follows from Proposition~\ref{prop:presheaf-exact} and
Corollary~\ref{cor:presheaf-sheaf}.
\pfend

\section{Low dimensional \v{C}ech cohomology}
\subsection{The group $\vH^0$}
Consider a simplicial space $\Mb$.
Let $\uUb$ be an open cover of $M_\com$, then,
using Notation~\ref{not:cocycle}, a 0-cocycle is
given by a family $(c_{\lambda_0})_{\lambda_0\in I_0}$,
with $c_{\lambda_0}\in \aA^0(U^0_{\lambda_0})$, and
\begin{equation}\label{eqn:0-cocycle}
0{\equal}(dc)_{\lambda_0\lambda_1\lambda_{01}}{\equal}
\teps_1^*c_{\lambda_1}-\teps_0^*c_{\lambda_0}
\end{equation}
on $U^1_\lambda{\equal}U^1_{\lambda_{01}}\cap\teps_0^{-1}(U^0_{\lambda_0})
\cap\teps_1^{-1}(U^0_{\lambda_1})$. Therefore, $\teps_1^*c_{\lambda_1}
{\equal}\teps_0^*c_{\lambda_0}$ on $\teps_0^{-1}(U^0_{\lambda_0})
\cap\teps_1^{-1}(U^0_{\lambda_1})$ for all $\lambda_0$,
$\lambda_1\in I_0$. Applying $\teta_0^*$ to both sides, we find
that $c_{\lambda_0}{\equal}c_{\lambda_1}$ on $U^0_{\lambda_0}\cap
U^0_{\lambda_1}$. Since $\aA^0$ is a sheaf, there exists a global
section $\varphi\in \aA^0(M_0)$ such that $c_{\lambda_0}$ is the
restriction of $\varphi$ to $U^0_{\lambda_0}$ for all $\lambda_0
\in I_0$. Now, Equation (\ref{eqn:0-cocycle})
is equivalent to $\teps_1^*\varphi
{\equal}\teps_0^*\varphi$. We have thus proved:

\begin{prop}\label{prop:H0}
Let $\aAb$ be an abelian sheaf on a simplicial space $\Mb$ and let $\uUb$
be an open cover of $\Mb$. Then
$$\vH^0(\Mb;\aAb){\equal}H^0(\uUb;\aAb){\equal}\Gamma_{\mathrm{inv}}(\aAb)
:{\equal}{\mathrm{Ker}}\,(\aA^0(M_0)\rightrightarrows\aA^1(M_1)).$$
\end{prop}

(Of course, in the case of a groupoid and an abelian $G$-sheaf,
a section in $\aA^0(G_0)$
is in $\Gamma_{\mathrm{inv}}(\aAb)$ if and only if it is an invariant
section in the usual sense, i.e. under the action of $G$.)

\subsection{The group $\vH^1$}
Consider a groupoid $G$. The cocycle relation in degree 1 is
$$\teps_0^*c_{\lambda_1\lambda_2\lambda_{12}}
-\teps_1^*c_{\lambda_0\lambda_2\lambda_{02}}+\teps_2^*
c_{\lambda_0\lambda_1\lambda_{01}}{\equal}0$$
on $U^2_{\lambda_0\lambda_1\lambda_2\lambda_{01}\lambda_{02}\lambda_{12}
\lambda_{012}}$.
Exactly the same method as in the preceding paragraph shows
that $c_{\lambda_0\lambda_1\lambda_{01}}$ does not depend on the choice
of $\lambda_{01}$, hence there exists a section $\varphi_{\lambda_0\lambda_1}
\in \aA^1(\teps_0^{-1}(U^0_{\lambda_0})\cap\teps_1^{-1}(U^1_{\lambda_1}))$
such that $c_{\lambda_0\lambda_1\lambda_{01}}$ is the restriction
to $U^1_{\lambda_0\lambda_1\lambda_{01}}$ of $\varphi_{\lambda_0\lambda_1}$.
The cocycle relation becomes
\begin{equation}\label{eqn:1-cocycle}
\teps^*_0\varphi_{\lambda_1\lambda_2}
-\teps^*_1\varphi_{\lambda_0\lambda_2}
+\teps^*_2\varphi_{\lambda_0\lambda_1}{\equal}0.
\end{equation}

Coboundaries are cocycles of the form
$\varphi_{\lambda_0\lambda_1}{\equal}\teps_0^*\alpha_{\lambda_1}
-\teps_1^*\alpha_{\lambda_0}$.

\begin{prop}
Let $G$ be a topological groupoid and $A$ be a $G$-module.
Denote by $\aAb$ the associated sheaf on $G_\com$.
Then $H^1(G_\com;\aAb)$ is the group of $G$-equivariant
locally trivial $A$-principal bundles over $G_0$.
\end{prop}

\pf
This is well-known. Let us give the proof for completeness.
Suppose we are given such a principal bundle $E\to G_0$. Choose
an open cover $\uU_0{\equal}(U^0_i)_{i\in I_0}$ of $G_0$ such that for all
$i$ there exists a (not necessarily equivariant) section
$x\mapsto \sigma_i(x)$ of the bundle $E\to G_0$ over the open subset
$U^0_i$. Let us define
$\varphi_{\lambda_0\lambda_1}(g)\in A_{r(g)}$ (for all
$g\in \teps^{-1}_0(U^0_{\lambda_0})
\cap \teps^{-1}_1(U^0_{\lambda_1})$) by
\begin{equation}\label{eqn:1-cocycle2}
\sigma_{\lambda_0}(r(g)) {\equal} g\sigma_{\lambda_1}(s(g))
+\varphi_{\lambda_0\lambda_1}(g).
\end{equation}
(We have denoted by ``$+$'' the action of $A$ on the $A$-principal
bundle $E$.)
If $(g,h)\in G_2$, then we get
\begin{eqnarray}\label{eqn:1-cocycle3}
\sigma_{\lambda_1}(r(h)) &{\equal}& h\sigma_{\lambda_2}(s(h))
+\varphi_{\lambda_1\lambda_2}(h)\\
\label{eqn:1-cocycle4}
\sigma_{\lambda_0}(r(gh)) &{\equal}& gh\sigma_{\lambda_2}(s(gh))
+\varphi_{\lambda_0\lambda_2}(gh).
\end{eqnarray}
Substitute (\ref{eqn:1-cocycle3}) in (\ref{eqn:1-cocycle2})
and compare to (\ref{eqn:1-cocycle4}) to get the cocycle relation
\begin{equation}\label{eqn:1-cocycle5}
\varphi_{\lambda_0\lambda_2}(gh){\equal}g\varphi_{\lambda_1\lambda_2}(h)
+\varphi_{\lambda_0\lambda_1}(g)
\end{equation}
(which is exactly (\ref{eqn:1-cocycle})).

If $\sigma'_i$ is another section, let $\alpha_i(x)
{\equal}\sigma_i(x)-\sigma'_i(x)$, then from (\ref{eqn:1-cocycle2})
and its analogue for $\sigma'$, we get by substraction
$\alpha_{\lambda_0}(r(g)){\equal}g\alpha_{\lambda_1}(s(g))
+\varphi_{\lambda_0\lambda_1}(g)-\varphi'_{\lambda_0\lambda_1}(g)$,
i.e. $\varphi'-\varphi{\equal}d\alpha$.

The above shows that any $G$-equivariant $A$-principal bundle which is
trivial over each open set $U^0_i$ defines an element of
$H^1(\uUb;\aAb)$. Moreover, passing to a finer cover obviously
doesn't change the \v{C}ech cohomology class, hence
any bundle as above defines a cohomology class in $\vH^1(G_\com;\aAb)$.
\par\medskip

Conversely, suppose we are given $(\varphi_{\lambda_0\lambda_1})$
satisfying the cocycle relation (\ref{eqn:1-cocycle5}). Define
$$E{\equal}(\coprod_{i\in I_0} U^0_i\times_{G_0} A)/\sim$$
with the identifications
$(\lambda_0,x,a)\sim (\lambda_1,x,a+\varphi_{\lambda_0\lambda_1}(x))$
and action
$$g\cdot[(\lambda_1,s(g),a)]{\equal}[(\lambda_0,r(g),
a-\varphi_{\lambda_0\lambda_1}(g))].$$
It is elementary to check that $E$ is a locally trivial $G$-equivariant
$A$-principal bundle, and that the associated 1-cocycle is indeed
$\varphi$.
\pfend

\subsection{The group $\vH^2$, extensions and the Brauer group}
Let $G$ be a topological groupoid, and let $A$ be a $G$-module.
Let us denote by $\mathrm{ext}(G,A)$ the set of extensions
of the form
$$A\stackrel{i}{\hookrightarrow} E \stackrel{\pi}{\twoheadrightarrow} G$$
such that the unit spaces of the groupoids $A$, $E$ and $G$ are all
equal to $G_0$, the maps $i$ and $\pi$ are the identity map on
$G_0$, and such that for all $\gamma\in E$ and all $a\in A_{s(\gamma)}$
we have
$$\gamma a\gamma^{-1}{\equal}\pi(\gamma)\cdot a.$$
In $\pi(\gamma)\cdot a$, the dot denotes the action of $G$ on the
$G$-module $A$.

Two such extensions $A\to E\to G$ and $A\to E'\to G$ are considered
equivalent if there is commutative diagram
$$\xymatrix{
A\ar[r]\ar[d]^{\text{Id}} &E\ar[r]\ar[d]& G\ar[d]^{\text{Id}}\\
A\ar[r] & E'\ar[r] & G
}$$
such that the map $E\to E'$ is a groupoid isomorphism.

There is a canonical extension in $\mathrm{ext}(G,A)$: let $E{\equal}A\times_{p,r} G
{\equal}\{(a,g)\in A\times G\vert\; p(a){\equal}r(g)\}$.
The source and range maps
in $E$ are $s(a,g){\equal}s(g)$, $r(a,g){\equal}r(g)$. The product is
$(a,g)(b,h){\equal}(a+g\cdot b,gh)$ (defined whenever $s(g){\equal}r(g)$;
the product in $A$ is written additively). The inclusion
$A\hookrightarrow E$ is $i(a){\equal}(a,p(a))$ and the projection
$\pi\colon E\to G$ is $\pi(a,g){\equal}g$.

Let us call this extension the
\emph{strictly trivial} extension.

\begin{prop}\label{prop:strictly-trivial}
Let $A\stackrel{i}{\hookrightarrow} E \stackrel{\pi}{\twoheadrightarrow} G$
be an element of $\mathrm{ext}(G,A)$.
The following are equivalent:
\begin{itemize}
\item[(i)] the extension is strictly trivial;
\item[(ii)] there exists a groupoid morphism $\sigma\colon G\to E$
which is a section of $\pi$;
\item[(iii)] there exists $\varphi\colon E\to A$
such that $\varphi(a\gamma){\equal}a\varphi(\gamma)$ for all
$(a,\gamma)\in A\times_{p,r}E$ and $\varphi(\gamma_1\gamma_2){\equal}
\varphi(\gamma_1)+\pi(\gamma_1)\cdot\varphi(\gamma_2)$ for all
composable pairs $(\gamma_1,\gamma_2)\in E^2$.
\end{itemize}
\end{prop}

\pf
Let us sketch the easy proof.

(i) $\implies$ (ii): take $\sigma(g){\equal}(r(g),g)$.

For (ii) $\implies$ (iii), let
$\varphi(\gamma){\equal}\gamma[\sigma\circ\pi(\gamma)]^{-1}$.

To prove (iii) $\implies$ (i), the map $\gamma\mapsto
(\varphi(\gamma),\pi(\gamma))$ is a groupoid isomophism
from $E$ onto $A\times_{p,r} G$.
\pfend

The set $\mathrm{ext}(G,A)$ is an abelian group. The (``Baer'')
sum $E_1\oplus E_2$ of
two extensions $A\to E_i\to G$ is given by the extension
$A\to E\to G$ with
$$E{\equal}\{(\gamma_1,\gamma_2)\in E_1\times E_2\vert\;
\pi_1(\gamma_1){\equal}\pi_2(\gamma_2)\}/\sim$$
where $(a\gamma_1,\gamma_2)\sim (\gamma_1,a\gamma_2)$
if $\pi_1(\gamma){\equal}\pi_2(\gamma_2)$ and $p(a){\equal}r(\gamma_1)$.
The map $\pi\colon E\to G$ is given by $\pi(\gamma_1,\gamma_2){\equal}
\pi_1(\gamma_1){\equal}\pi_2(\gamma_2)$, and the inclusion $i\colon A\to E$
is
$$i(a){\equal}(i_1(a),p(a))\sim(p(a),i_2(a)).$$
Finally, the inverse of the extension is
$$A\stackrel{i'}{\hookrightarrow} \bar{E}
\stackrel{\pi}{\twoheadrightarrow} G$$
where $\bar{E}{\equal}E$ as a groupoid, but $i'(a){\equal}i(-a)$. Denoting by
$\bar{\gamma}$ the element in $\bar{E}$ which is the same as the
element $\gamma\in E$, this means that
$\overline{a\gamma}{\equal}(-a)\bar{\gamma}$.
To check that
$E\oplus \bar{E}$ is strictly trivial, just note that for all
$g\in G$, the element $(\gamma,\bar{\gamma})\in E\oplus \bar{E}$, does
not depend on the choice of $\gamma\in E$
such that $\pi(\gamma){\equal}g$, since
$(a\gamma,\overline{a\gamma}){\equal}(a\gamma,(-a)\bar{\gamma})\sim(\gamma,
\bar{\gamma})$. Therefore, $g\mapsto \sigma(g){\equal}(\gamma,\bar{\gamma})$
defines a cross-section of $\pi$.

\begin{example}
When $A{\equal}G_0\times \T$, and $G$ does not act on $\T$, we obtain the group
${\mathrm{Tw}}(G)$ of twists of $G$ \cite{KMRW98}.
\end{example}

It is clear that $\mathrm{ext}(G,A)$ is covariant with respect to
$G$-module morphisms. It is also contravariant with respect to
groupoid morphisms. Indeed, let $f\colon G'\to G$ be a groupoid
morphism, and let $A'{\equal}f^*A{\equal}\{(x,a)\in G'_0\times A\vert\;
f(x){\equal}p(a)\}$. Then $A'$ is a $G'$-module with respect to the action
$g'\cdot (x,a){\equal}(r(g'),f(g)\cdot a)$, and there is a ``pull-back''
morphism
$$f^*\colon \mathrm{ext}(G,A)\to \mathrm{ext}(G',f^*A)$$
defined as follows. Let $A\to E\to G$ be an element of
$\mathrm{ext}(G,A)$, then its pull-back by $f$ is the extension
$$A\stackrel{i'}{\hookrightarrow} E' \stackrel{\pi'}{\twoheadrightarrow} G'$$
where $E'{\equal}\{(\gamma,g')\in E\times G'\vert\;
\pi(\gamma){\equal}f(g')\}$, $\pi'(\gamma,g'){\equal}g'$,
$i'(a){\equal}(i(a),p(a))$. The groupoid structure on $E'$ is the one
induced from the product groupoid $E\times G'$.

More generally, suppose that $G'_0\stackrel{\rho}{\leftarrow}
Z\stackrel{\sigma}{\to}G_0$ is a generalized morphism from
$G'$ to $G$ (see \ref{sec:Morita}). Put $A'{\equal}Z\times_G A
:{\equal}\{(z,a)\in Z\times A\vert\; \sigma(z){\equal}p(a)\}/\sim$,
where $(zg,g^{-1}a)\sim (z,a)$ for all triples $(z,a,g)\in Z\times A\times G$
such that
$\sigma(a){\equal}p(a){\equal}r(g)$. It is obvious that $A'$ is a $G'$-module
with sum $(z,a)+(z,b){\equal}(z,a+b)$ and left $G'$-action
$g'(z,a){\equal}(g'z,a)$.

The slight defect of the group $\mathrm{ext}(G,A)$ is that it is
not invariant by Morita equivalence. To remedy this, let us define

\begin{defi}\label{def:Ext}
$$\mathrm{Ext}(G,A){\equal}\lim_{\uU}\mathrm{ext}(G[\uU],A[\uU])$$
where $\uU$ runs over open covers of $G_0$ (see notation (\ref{eqn:GU})).
\end{defi}

By construction, the group $\mathrm{Ext}(G,A)$ is invariant under Morita
equivalence (see Proposition~\ref{prop:Morita}).
\par\medskip

Let us now come to the relation between 2-cohomology and extensions:

\begin{prop}\label{prop:2-cohomology}
Let $G$ be a topological groupoid, $A$ a $G$-module, and $\aAb$ the
associated sheaf over $G_\com$.
\begin{itemize}
\item[(a)] For each open cover $\uUb$ of
$G_\com$, there is a canonical group isomorphism
\begin{equation}\label{eqn:extension-cover}
\mathrm{ext}_\uU(G[\uU_0],A[\uU_0])\cong H^2(\uUb;\aAb),
\end{equation}
where $\mathrm{ext}_\uU(G[\uU_0],A[\uU_0])$ denotes the subgroup of
elements of $\mathrm{ext}(G[\uU_0],A[\uU_0])$ consisting of extensions
$A[\uU_0]\to E\stackrel{\pi}{\to} G[\uU_0]$ such that $\pi$ admits
a continuous lifting over each open set $U^1_\lambda\subset G_1$
($\lambda\in \Lambda_1$).
\item[(b)] (\ref{eqn:extension-cover})
induces an isomorphism
$$\mathrm{Ext}(G,A)\cong \vH^2(\uUb;\aAb).$$
\end{itemize}
\end{prop}

\pf
As in the previous subsection, one easily sees that a 2-cocycle
in $Z^2(\uUb;\aAb)$ is given by a family
$$\varphi{\equal}(\varphi_{\lambda_0\lambda_1\lambda_2\lambda_{01}
\lambda_{02}\lambda_{12}})$$
such that each term is a continuous function $(g,h)\mapsto
\varphi_\lambda(g,h)\in A_{r(g)}$, defined on the set of pairs
$(g,h)$ such that $r(g)\in U^0_{\lambda_0}$, $s(g)\in
U^0_{\lambda_1}$, $s(h)\in U^0_{\lambda_2}$,
$g\in U^1_{\lambda_{01}}$, $gh\in U^1_{\lambda_{02}}$,
$h\in U^1_{\lambda_{12}}$. The $\varphi$'s satisfy the cocycle identity
\begin{eqnarray}\label{eqn:2-cocycle}
\lefteqn{g\varphi_{\lambda_1\lambda_2\lambda_3\lambda_{12}
\lambda_{13}\lambda_{23}}(h,k)
-\varphi_{\lambda_0\lambda_2\lambda_3\lambda_{02}
\lambda_{03}\lambda_{23}}(gh,k)}\\\nonumber
&&+\varphi_{\lambda_0\lambda_1\lambda_3\lambda_{01}
\lambda_{03}\lambda_{13}}(g,hk)
-\varphi_{\lambda_0\lambda_1\lambda_2\lambda_{01}
\lambda_{02}\lambda_{12}}(g,h){\equal}0
\end{eqnarray}

Let us consider a cover $\vVb$ of $G[\uU]$, $\vV_n{\equal}(V^n_j)_{j\in J_n}$,
such that
\begin{itemize}
\item $J_0{\equal}\{\mathrm{pt}\}$ and $\vV_0$ is the cover consisting of the
unique open set $\coprod_{i\in I_0}U^0_i$;
\item $J_1{\equal}I_0\times I_0\times I_1$ with
$V^1_{ijk}{\equal}\{(i,g,j)\vert\; r(g)\in U^0_i,\;s(g)\in U^0_j,\;
g\in U^1_k\}$;
\item $V^n$ arbitrary $\forall n\ge 2$.
\end{itemize}

Consider the group $Z^2(\vVb;{\aA}^{'\com})$, where ${\aA}^{'\com}$
is the pull-back of the sheaf $\aAb$ by $G[\uU]\to G$.
As above, it consists of families $\psi_{\mu_{01}\mu_{02}\mu_{12}}$
satisfying the cocycle identity
\begin{eqnarray}\label{eqn:2-cocycle2}
\lefteqn{g\psi_{\mu_{12}
\mu_{13}\mu_{23}}(h,k)
-\psi_{\mu_{02}
\mu_{03}\mu_{23}}(gh,k)}\\\nonumber
&&+\psi_{\mu_{01}
\mu_{03}\mu_{13}}(g,hk)
-\psi_{\mu_{01}
\mu_{02}\mu_{12}}(g,h){\equal}0
\end{eqnarray}

We show that $Z^2(\uUb;\aAb)\cong Z^2(\vVb;{\aA}^{'\com})$, where ${\aA}^{'\com}$
is the pull-back of the sheaf $\aAb$ by $G[\uU]\to G$.

In one direction, let $\psi\in Z^2(\vVb;{\aA}^{'\com})$. For all
$\lambda{\equal}
(\lambda_0,\lambda_1,\lambda_2,\lambda_{01},\lambda_{02},\lambda_{12})$
in $I_0^3\times I_1^3$, define
\begin{eqnarray}\label{eqn:lambda-mu}
\mu_{01}&{\equal}&(\lambda_0,\lambda_1,\lambda_{01})\\\nonumber
\mu_{02}&{\equal}&(\lambda_0,\lambda_2,\lambda_{02})\\\nonumber
\mu_{12}&{\equal}&(\lambda_1,\lambda_2,\lambda_{12})
\end{eqnarray}
and $\varphi_\lambda{\equal}\psi_{\mu_{01}\mu_{02}\mu_{12}}$.

In the other direction, if we are given a 2-cocycle $\varphi
\in Z^2(\uUb;\aAb)$, we want to define a 2-cocycle
$\psi\in Z^2(\vVb;{\aA}^{'\com})$. Given $\mu{\equal}(\mu_{01},\mu_{02},
\mu_{12})\in J_1$, write $\mu_{ab}{\equal}(i_{ab},j_{ab},k_{ab})$.
Then $V^1_\mu\ne\emptyset$ implies that $j_{01}{\equal}i_{12}$,
$i_{01}{\equal}i_{02}$, $j_{02}{\equal}j_{12}$, hence there exists
$$\lambda{\equal}(\lambda_0,\lambda_1,\lambda_2,\lambda_{01},
\lambda_{02},\lambda_{12})\in I_0^3\times I_1^3$$
such that (\ref{eqn:lambda-mu}) holds.
We then define $\psi_\mu{\equal}\varphi_\lambda$.

Comparing (\ref{eqn:2-cocycle}) and (\ref{eqn:2-cocycle2}),
we see that $Z^2(\uUb;\aAb)\cong Z^2(\vVb;{\aA}^{'\com})$.
Moreover, it is not hard to check that this
induces an isomorphism
\begin{equation}\label{eqn:cohomology-cover}
H^2(\uUb;\aAb)\cong H^2(\vVb;{\aA}^{'\com}).
\end{equation}
\par\medskip

To prove the first part of the proposition, then, we can (after passing
to the groupoid $G[\uU_0]$), suppose that $\uU_0$ consists of the
unique open set $G_0$.

Consider an extension in $\mathrm{ext}_{\uU}(G,A)$
$$A\hookrightarrow E \twoheadrightarrow G.$$
For each $i\in I_1$, consider a continuous section $\sigma_i
\colon U^1_i\to E$. Define a cochain $\varphi$ by the
equation
\begin{equation}\label{eqn:extension-2-cocycle}
\sigma_{\lambda_{01}}(g)
\sigma_{\lambda_{12}}(h)
{\equal}\varphi_{\lambda_{01}\lambda_{02}\lambda_{12}}(g,h)
\sigma_{\lambda_{02}}(gh).
\end{equation}
To see that it is indeed a cocycle, just write
$$(\sigma_{\lambda_{01}}(g)\sigma_{\lambda_{12}}(h))
\sigma_{\lambda_{23}}(k){\equal}\sigma_{\lambda_{01}}(g)(\sigma_{\lambda_{12}}(h)
\sigma_{\lambda_{23}}(k))$$
and substitute relations like (\ref{eqn:extension-2-cocycle}) to obtain

\begin{eqnarray}\label{eqn:2-cocycle3}
\lefteqn{g\varphi_{\lambda_{12}
\lambda_{13}\lambda_{23}}(h,k)
-\varphi_{\lambda_{02}
\lambda_{03}\lambda_{23}}(gh,k)}\\\nonumber
&&+\varphi_{\lambda_{01}
\lambda_{03}\lambda_{13}}(g,hk)
-\varphi_{\lambda_{01}
\lambda_{02}\lambda_{12}}(g,h){\equal}0.
\end{eqnarray}

Suppose that $\sigma'_i$ is another continuous lifting and let
$\alpha_i\colon U^1_i\to A$ such that
\begin{equation}\label{eqn:alphai}
\sigma'_i(g){\equal}\alpha_i(g)\sigma_i(g).
\end{equation}
Define $\varphi'$ by
\begin{equation}\label{eqn:extension-2-cocycle2}
\sigma'_{\lambda_{01}}(g)
\sigma'_{\lambda_{12}}(h)
{\equal}\varphi'_{\lambda_{01}\lambda_{02}\lambda_{12}}(g,h)
\sigma'_{\lambda_{02}}(gh).
\end{equation}
Substituting (\ref{eqn:alphai}) in (\ref{eqn:extension-2-cocycle2})
and comparing with (\ref{eqn:extension-2-cocycle}), we find
$$(\varphi'-\varphi)_{\lambda_{01}\lambda_{02}\lambda_{12}}(gh)
{\equal}g\alpha_{\lambda_{12}}(h)-\alpha_{\lambda_{02}}(gh)+
\alpha_{\lambda_{01}}(g),$$
i.e. $\varphi'-\varphi{\equal}d\alpha$. This proves that an extension
in $\mathrm{ext}_{\uU}(G,A)$ determines a unique
cohomology class in $H^2(\uUb;\aAb)$.
\par\medskip
Conversely, given a cocycle $\varphi_{\lambda_{01}\lambda_{02}\lambda_{12}}$,
we want to construct an extension
$$A\to E\to G.$$
The idea is to set
\begin{equation}\label{eqn:cocycle-extension}
E{\equal}\coprod_{i\in I_1}\{(a,g,i)|\; a\in A,\;g\in U^1_i,\;p(a){\equal}r(g)\}/\sim,
\end{equation}
with the product law
\begin{equation}\label{eqn:product-law}
[a,g,\lambda_{01}][b,g,\lambda_{12}]{\equal}[a+g\cdot b+
\varphi_{\lambda_{01}\lambda_{02}\lambda_{12}}(g,h),gh,\lambda_{02}].
\end{equation}
To determine the correct equivalence relation in (\ref{eqn:cocycle-extension}),
we note that if $[a,x,i]$ represents a unit element in the groupoid $E$,
then from the product law (\ref{eqn:product-law}) we necessarily
have $[a,x,i]{\equal}[a,x,i][a,x,i]{\equal}[2a+\varphi_{iii}(x,x),x,i]$, thus
$[-\varphi_{iii}(x,x),x,i]$ must be the unit element.
Using again (\ref{eqn:product-law}), we necessarily have
$$[-\varphi_{iii}(r(g),r(g)),r(g),i][a,g,k]{\equal}
[-\varphi_{iii}(r(g),r(g))+a+\varphi_{ijk}(r(g),g),g,j]$$
thus we necessarily have
\begin{equation}\label{eqn:equivalence-relation}
(a,g,k)\sim(-\varphi_{iii}(r(g),r(g))+a+\varphi_{ijk}(r(g),g),g,j)
\end{equation}
Conversely, we want to show that (\ref{eqn:equivalence-relation})
defines an equivalence relation. We claim that $\psi_{kj}(g)
{\equal}-\varphi_{iii}(r(g),r(g))+\varphi_{ijk}(r(g),g)$ does not depend
of the choice of $i$.

Let us denote $x{\equal}r(g)$. Apply (\ref{eqn:2-cocycle3}) to
$(x,x,g)$ instead of $(g,h,k)$:
\begin{equation}\label{eqn:2-cocyclexxg}
\varphi_{ijk}(x,g)-\varphi_{\ell m k}(x,g)
+\varphi_{nmj}(x,g)-\varphi_{n\ell i}(x,x){\equal}0.
\end{equation}
Taking $g{\equal}x$ and $j{\equal}k{\equal}\ell{\equal}m{\equal}n$ we find
\begin{equation}\label{eqn:phi-xxx}
\varphi_{ijj}(x,x){\equal}\varphi_{jji}(x,x){\equal}\varphi_{iii}(x,x).
\end{equation}
Taking $k{\equal}j$ and $n{\equal}\ell$ in (\ref{eqn:2-cocyclexxg}) and
using (\ref{eqn:phi-xxx}), we get
\begin{equation}\label{eqn:phi-xg}
\varphi_{ijj}(x,g){\equal}\varphi_{\ell \ell i}(x,x){\equal}\varphi_{iii}(x,x).
\end{equation}
Then, take $m{\equal}j$ and $n{\equal}\ell$ in (\ref{eqn:2-cocyclexxg}) and
use (\ref{eqn:phi-xxx}) and (\ref{eqn:phi-xg}):
\begin{eqnarray*}
\varphi_{ijk}(x,g)-\varphi_{\ell jk}(x,g)+
\varphi_{\ell j j}(x,g)-\varphi_{\ell\ell i}(x,x) {\equal} 0\\
\varphi_{ijk}(x,g)-\varphi_{\ell jk}(x,g)+\varphi_{\ell\ell\ell}(x,x)
-\varphi_{iii}(x,x){\equal}0.
\end{eqnarray*}
This proves our claim that $\psi_{kj}$ is well-defined. Moreover, taking
$n{\equal}\ell{\equal}i$ in (\ref{eqn:2-cocyclexxg}) we get
\begin{equation}\label{eqn:psi}
\psi_{jk}(g)-\psi_{mk}(g)+\psi_{mj}(g){\equal}0.
\end{equation}
It follows that
\begin{eqnarray*}
\psi_{jj}&{\equal}&0\quad(\mbox{use }(\ref{eqn:psi})\mbox{ for }k{\equal}j)\\
\psi_{kj}&{\equal}&-\psi_{jk}\quad(\mbox{use }(\ref{eqn:psi})\mbox{ for }m{\equal}k)\\
\psi_{jm}&{\equal}&\psi_{jk}+\psi_{km}.
\end{eqnarray*}
Therefore, (\ref{eqn:equivalence-relation})
defines an equivalence relation.

It is then elementary to check that
(\ref{eqn:product-law}) endows $E$ with a groupoid structure
such that the obvious extension
$$A\to E\stackrel{\pi}{\to} G$$
is an element of $\mathrm{ext}(G,A)$, and that
$\pi$ admits a continuous lifting $\sigma_i\colon U^1_i\to G$
defined by $\sigma_i(g){\equal}[0,g,i]$; and that the associated cocycle
is precisely $\varphi$. We leave these easy verifications to the reader.

\par\bigskip
To prove the second part of the proposition, we first pass to the
inductive limit over all open covers $\uU_1$ of $G_1$ (leaving
$\uU_0$ fixed) to find that
$$\mathrm{ext}(G[\uU_0],A[\uU_0]){\equal}\lim_{\uU_1}H^2(\uUb;\aAb)$$
and then take the inductive limit over $\uU_0$.
\pfend

\begin{rem}\label{rem:cohomology-cover}
By the same method we used to show (\ref{eqn:cohomology-cover}),
one can show that for each open cover $\uUb$ of $G_\com$, and each
sheaf $\aAb$ over $G_\bullet$, the canonical morphism
$f\colon G[\uU_0]\to G$ induces an isomorphism
$$H^n(\uUb;\aAb)\cong H^n(\vVb;{\aA'}^\com)$$
where $\vV_0$ is the cover consisting of the unique open set
$\coprod_{i\in I}U^0_i$ and the open cover $\vV_n{\equal}
(V^n_{i_0,\ldots,i_n,j})_{i_0,\ldots,i_n,j\in I_0^{n+1}\times I_n}$
of $G[\uU]_n$ is defined by
\begin{eqnarray*}
V^n_{i_0,\ldots,i_n,j}&{\equal}&
\{(i_0,\ldots,i_n,g_1,\ldots,g_n)\vert\;
r(g_1)\in U^0_{i_0},\;s(g_1)\in U^0_{i_1},\ldots,\\
&&\qquad s(g_n)\in U^0_{i_n},\;
(g_1,\ldots,g_n)\in U^n_j\}
\end{eqnarray*}
(see Remark~\ref{rem:GU}).
\end{rem}

\begin{rem}
For completeness, it would remain to examine the relation
between the \emph{sheaf} 2-cohomology groups and extensions
of non-paracompact groupoids, since it is not obvious whether
$H^2$ is equal or not to $\vH^2$ in this case. However, we won't
develop this, due to lack of applications.
\end{rem}

\begin{coro}
If $G$ is a locally compact Hausdorff groupoid with Haar system
then $\vH^2(G_\com;\T)
\cong {\mathrm{Ext}}(G,G_0\times \T)$ is the Brauer group of $G$
\end{coro}

\pf
Use for instance Proposition~2.13
and remarks preceding Proposition~2.29 of \cite{TXL03}.
\pfend

\section{Comparison with Moore's cohomology}
Recall \cite{Moo76}
that if $G$ is a locally compact group and $A$ is a Polish
abelian group (i.e., as a topological space, $A$ admits a separable
complete metric), then $A$ is a $G$-module if $G$ acts (continuously)
by automorphisms on $A$.

Given a Polish $G$-module $A$, let $I(A)$ be the set of
$\mu$-measurable functions from $G$ to $A$ ($\mu$ being the Haar measure),
modulo equality almost everywhere. Then $I(A)$ is again a Polish
$G$-module, with action $(\gamma\cdot f)(x){\equal}\gamma f(\gamma^{-1}x)$
(caution: our definition is different but isomorphic to Moore's
definition of $I(A)$).

The $G$-module $A$ embeds in $I(A)$ via the obvious map
$$i_A\colon A\to I(A),\quad (i_A(a))(x){\equal}a.$$
Let $U(A){\equal}I(A)/A$. Then, using measurable cocycles, Moore defined
cohomology groups $H^n(G,A)$ which are
are characterized by the proposition below, where $I(A)$ is
defined as above and $F(A){\equal}A^G$ (the sub-module of $G$-fixed points).

\begin{prop}\label{prop:derived-functors}
Let $\cC_1$ and $\cC_2$ be two abelian categories. Suppose that $F$ is
a left-exact functor from $\cC_1$ to $\cC_2$, that $I\colon \cC_1\to \cC_2$
is a functor and $i_A\colon A\hookrightarrow I(A)$ is a natural injection.
Then
\begin{itemize}
\item[(a)] there exists, up to isomorphism, at most one sequence of functors
$H^n\colon \cC_1\to \cC_2$ such that
\begin{itemize}
\item[1)] $H^0{\equal}F$
\item[2)] Any exact sequence $0\to A'\to A\to A''\to 0$ induces a
natural long exact sequence $0\to H^0(A')\to H^0(A)\to H^0(A'')
\stackrel{\partial}{\to} H^1(A')\to \ldots$
\item[3)] $H^n(I(A)){\equal}0$ for all $A$ and for all $n\ge 1$.
\end{itemize}
\item[(b)] If $I$ is an exact functor and $I(i_A){\equal}i_{I(A)}$ for all
$A$, then 3) may be replaced by
\ \ 3') $H^1(I(A)){\equal}0$ $\forall A$.
\item[(c)] If moreover $F(I(A))\to F(U(I(A)))$ is surjective for all
$A$, where $U(A){\equal}I(A)/A$, then there exists a sequence of functors satisfying
1), 2) and 3).
\end{itemize}
\end{prop}

\pf
a) This is essentially \cite[Theorem~1]{Ati-Wal67} or
\cite[Theorem~2]{Moo76}: using the long exact sequence associated to
$0\to A\stackrel{i_A}{\to} I(A)\to U(A)\to 0$, one gets
$H^n(A)\cong H^{n-1}(U(A))$ for $n\ge 2$ and $H^1(A){\equal}
\mathrm{coker}\,(F(I(A))\to F(U(A)))$, thus $H^n$ is uniquely
determined by induction on $n$.

b) If $I$ is exact then
$$0\to I(A)\stackrel{I(i_A)}{\to} I(I(A))\to I(U(A))\to 0$$
is exact, and
$$0\to I(A)\stackrel{i_{I(A)}}{\to} I(I(A))\to U(I(A))\to 0$$
is exact by definition. Therefore, the assumption $I(i_A){\equal}i_{I(A)}$
implies $I(U(A)){\equal}U(I(A))$ canonically. Thus, 3')
implies that for $n\ge 2$,
$$H^n(I(A)){\equal}H^1(U^{n-1}I(A)){\equal}H^1(IU^{n-1}(A)){\equal}0.$$

c) Define a resolution $0\to A\to A_0\to A_1\to\cdots$ by
$A_0{\equal}I(A$), and $A_{n+1}{\equal}I(A_n/A_{n-1})$. The map $A_n\to A_{n+1}$
is the composition $A_n\to A_n/A_{n-1}\stackrel{i}{\to}
I(A_n/A_{n-1})$. Define $H^n(A)$ to be the cohomology of the
complex $F(A_0)\to F(A_1)\to\cdots$ and let us check that properties
1), 2) and 3') hold.

1) $H^0(A){\equal}\mathrm{Ker}\,(F(A)\to F(I(I(A)/A)))
{\equal}\mathrm{Ker}\,(F(I(A))\to F(I(A)/A))$ since $F$ preserves
injectivity of morphisms, and $I(A)/A\to I(I(A)/A)$ is
injective. Using left exactness of $F$,
we see that $H^0(A){\equal}F(\mathrm{Ker}(I(A)\to I(A)/A)){\equal}F(A))$.

2) Since $I$ is an exact functor, we have an exact sequence of complexes
$0\to A'_*\to A_*\to A''_*\to 0$, hence the conclusion by the Snake
lemma.

3') $H^1(I(A)){\equal}F(I(I(A)))/F(U(I(A))){\equal}0$.
\pfend

For instance, in part a) of the proposition, if $I(A)$ is an injective
object for all $A$ then $H^n$ are the right derived functors of $F$.
\par\smallskip

We are now ready to prove:
\begin{prop}
Let $G$ be a locally compact group.
Let $A$ be a Polish $G$-module and let $\aAb$ be the associated sheaf
on $G_\com$ (see Definition~\ref{def:G-module} and below).
Then $\vH^*(G_\com;\aAb)\cong H^*(G,A)$.
\end{prop}

\pf
We just need to check that $\vH^*(G_\com;\cdot)$ satisfies the
conditions 1)--3) of Proposition~\ref{prop:derived-functors}.

1) $\vH^0(G_\com;\aAb){\equal}A^G$ was proved in
Proposition~\ref{prop:H0} and 2) is Corollary~\ref{cor:long-exact-cech}.

$I$ is an exact functor \cite[Proposition~9]{Moo76}, and it is obvious
that $I(i_A){\equal}i_{I(A)}$. It thus remains to show that $\vH^1(G_\com;
\bBb){\equal}0$ if $\bBb$ is the sheaf on $G_\com$ associated to
the $G$-module $I(A)$. Recall that a 1-cocycle $\varphi$
is a continuous function $\varphi\colon G\to I(A)$ satisfying
$$g_1\varphi(g_2)-\varphi(g_1g_2)+\varphi(g_1){\equal}0$$
(see (\ref{eqn:1-cocycle5})), hence $\varphi$ is a cocycle in the
Moore complex. But $H^1_{\text{Moore}}(G,I(A)){\equal}0$, hence there exists
$\psi\in I(A)$ such that $\varphi(g){\equal}g\psi-\psi$ for all
$g\in G$. Therefore, $\varphi{\equal}0$ in $H^1(\uUb;\bBb)$.
\pfend

\begin{rem}
One can easily construct explicitly the isomorphism
$\vH^n(G_\com;\aAb)\to H^n_{\text{Moore}}(G,A)$. Take
$\varphi\in Z^n(\uUb;\aAb)$ where $\uUb$ is an open cover. Choose
measurable maps $\theta_k\colon G^k\to I_k$ such that
$x\in U^k_{\theta_k(x)}$ for all $x\in G^k$, and define
$$c(g_1,\ldots,g_n){\equal}\varphi_{\lambda(g_1,\ldots,g_n)}
(g_1,\ldots,g_n)$$
where $\lambda(g_1,\ldots,g_n)\in \Lambda_n$ is defined by
$\lambda(g_1,\ldots,g_n)(f){\equal}\theta_r(\tf(g_1,\ldots,g_n)) \in I_k$
for every $f\in \hom_{\Delta'}(k,n)$.
\end{rem}

\begin{rem}
One might wonder if it is possible to define, for every locally
compact groupoid and every (say, locally compact) $G$-module $A$,
an analogue of the Moore complex, using measurable cochains
$c(g_1,\ldots,g_n)\in A_{r(g_1)}$. In order to get
the usual cohomology groups when $G$ is a space, one should
probably use the sheaf over $G_0$ of functions
$c(g_1,\ldots,g_n)\in A_{r(g_1)}$ which are measurable in the
``leaf'' direction and continuous in the ``transverse'' direction.
Since this approach doesn't seem simpler or more useful than
\v{C}ech or sheaf cohomology of simplicial spaces, we won't develop this
further.
\end{rem}

\section{Comparison with sheaf cohomology and Haefliger's cohomology}
Let $\Mb$ be a simplicial space and $\aAb$ an abelian sheaf on $\Mb$.
By definition \cite{Del74}, the cohomology groups $H^n(\Mb;\aAb)$
are the derived functors of the functor $\Gamma_{\text{inv}}
(\Mb;\aAb){\equal}\mathrm{Ker}\,(\aA^0(M_0)\rightrightarrows\aA^1(M_1))$,
thus they coincide with Haefliger's cohomology groups in the case
of \et\ groupoids \cite{Hae79}.

A practical way of calculating the cohomology groups is to take a
resolution $(\lL^{\com,q})_{q\in\N}$ of $\aAb$ such that
$H^n(M_p;\lL^{p,q}){\equal}0$ $\forall n\ge 1$, $\forall p,q\ge 0$ and
take the cohomology of the double complex
$(\lL^{p,q}(M_p))$, where the first differential is
$d'{\equal}\sum_{k{\equal}0}^{p+1}(-1)^k\teps_k^*$ and the second differential
$d''$ is the differential in the resolution
$$\aA^p\to \lL^{p,0}\to \lL^{p,1}\to\cdots$$
(See \cite[\S 5.2.3]{Del74} in the general case,
\cite[\S 2.7]{Cra-Moe} or \cite{Hae79} in the case
of \et\ groupoids).

In this section, we show:

\begin{prop}\label{prop:cech-sheaf}
Let $\Mb$ be a paracompact simplicial space, and $\aAb$ an abelian sheaf
on $\Mb$. Then $H^*(\Mb;\aAb)\cong \vH^*(\Mb;\aAb)$. In particular,
$\vH^*(G_\com;\aAb)$ are Haefliger's cohomology groups if
$G$ is an \et\ paracompact groupoid and $\aAb$ is an abelian $G$-sheaf.
\end{prop}

We will again use Proposition~\ref{prop:derived-functors}.
Consider $I(\aA)^\com$ the sheaf such that
$I(\aA)^n(U)$ consists of maps $f$ (continuous or not) from
$\teps_0^{-1}(U)$ to $\aA^{n+1}$ such that $f(x)\in \aA^{n+1}_x$
for all $x\in \teps_0^{-1}(U)$. We will need the

\begin{lem}\label{lem:I(A)-acyclic}
For any open cover $\uUb$ and all $n\ge 1$, $H^n(\uUb;I(\aA)^\com)
{\equal}\{0\}$.
\end{lem}

\pf
Let us first explain the simplicial structure on the sheaf $I(\aA)^\com$.
Given $f\in \hom_\Delta(k,n)$, $U\subset M_k$, $V\subset M_n$
such that $\tf(V)\subset U$ and a section $\varphi$ of $I(\aA)^k$ over
$U$, we have to produce a section $\tf\varphi\in \Gamma(V,I(\aA)^n)$.

Define $f'\in\hom_\Delta(k+1,n+1)$ such that
\begin{equation}\label{eqn:f'}
f'(0){\equal}0\mbox{ and }\eps_0\circ f
{\equal}f'\circ\eps_0\colon [k]\to [n+1].
\end{equation}
Since
$\tf(V)\subset U$ we have $\tf'(\teps_0^{-1}(V))\subset
\teps_0^{-1}(U)$, thus we get a section
$$x\in \teps_0^{-1}(V)\mapsto
\varphi(\tf'(x))\in \aA^{k+1}_{\tf'(x)}
\mapsto (\tf')^*\varphi(\tf'(x))\in \aA^{n+1}_x.$$
\par\medskip

Now, let us show that the complex $C^*(\uUb;I(\aA)^\com)$
is homotopically trivial. First, for all $n\in\N$ and all
$x\in M_n$, let us choose $\theta(x)\in I_n$ such that
$x\in U^n_{\theta(x)}$. We define a homotopy
$$H\colon C^n(\uUb;I(\aA)^\com)\to C^{n-1}(\uUb;I(\aA)^\com)$$
by $(H\varphi)_\lambda(x){\equal}\teta_0^*\varphi_{\lambda'_x}(\teta_0(x))$,
$\forall x\in\teps_0^{-1}(U_\lambda)\subset M_n$,
$\forall \lambda\in \Lambda_{n-1}$, and $\lambda'_x\in \Lambda_n$
is defined as follows: for all $f\in \hom_{\Delta'}(k,n)$, let
$$\lambda'_x(f){\equal}\left\{
\begin{array}{ll}
\lambda(\eta_0\circ f)&\mbox{if }f(0)\ne 0\\
\theta(\tf(x))&\mbox{if }f(0){\equal}0
\end{array}
\right\}\in I_k.$$
Let $\varphi\in C^n(\uUb;I(\aA)^\com)$. Let
us compute $(dH+Hd)\varphi$ and compare it with $\varphi$.
We have
$$(d\varphi)_\lambda(x){\equal}\sum_{k{\equal}0}^{n+1}
(-1)^k{\teps_k}^{'*}\varphi_{\teps_k(\lambda)}(\teps'_k(x))\quad
\in \aA^{n+2}_x.$$
The meaning of this formula is the following: we take $x\in M_{n+2}$,
then its image $\teps'_k(x)$ (see notation~(\ref{eqn:f'})) is in
$M_{n+1}$. Its image $\varphi_{\teps_k(\lambda)}(\teps'_k(x))$
belongs to $\aA^{n+1}_{\teps'_k(x)}$ is restricted (see (\ref{eqn:f*}))
by $\teps'_k\colon M_{n+2}\to M_{n+1}$ to an element of
$\aA^{n+2}_x$.
\par\smallskip
Actually, we have $\varepsilon'_k{\equal}\varepsilon_{k+1}$ since
$\varepsilon_0\circ\varepsilon_k {\equal} \varepsilon_{k+1}\circ
\varepsilon_0\colon [n]\to [n+2]$, hence
\begin{eqnarray*}
(d\varphi)_\lambda(x)&{\equal}&\sum_{k{\equal}0}^{n+1} (-1)^k\teps_{k+1}^*
\varphi_{\teps_k(\lambda)}(\teps_{k+1}(x))\\
(dH\varphi)_\lambda(x)&{\equal}&\sum_{k{\equal}0}^n (-1)^k\teps^*_{k+1}
(H\varphi)_{\teps_k(\lambda)}(\teps_{k+1}(x))\\
&{\equal}&\sum_{k{\equal}0}^n (-1)^k\teps^*_{k+1}\teta^*_0
\varphi_{\teps_k(\lambda)'_{\teps_{k+1}}(x)}(\teta_0\circ\teps_{k+1}(x))\\
(Hd\varphi)_\lambda(x)&{\equal}&\teta_0^*(d\varphi)_{\lambda'_x}(\teta_0(x))\\
&{\equal}&\sum_{k{\equal}0}^{n+1} (-1)^k\teta_0^*\teps_{k+1}^*
\varphi_{\teps_k(\lambda'_x)}(\teps_{k+1}\teta_0(x)).
\end{eqnarray*}
In the last sum, for $k{\equal}0$ we get $\varphi_{\teps_0(\lambda'_x)}(x)$.
Now, $(\teps_0(\lambda'_x))(f){\equal}\lambda'_x(\varepsilon_0\circ f)
{\equal}\lambda(\eta_0\circ\varepsilon_0\circ f){\equal}\lambda(f)$, thus
the term for $k{\equal}0$ is just $\varphi_\lambda(x)$.

To show that the other terms in the sum $dH+Hd$ cancel out, we just
need to check that for all $k\ge 1$,
\begin{itemize}
\item[(a)] $\eta_0\circ\varepsilon_{k+1}{\equal}\varepsilon_k\circ\eta_0$, and
\item[(b)] $\teps_k(\lambda'_x){\equal}\teps_{k-1}(\lambda)'_{\teps_k(x)}$.
\end{itemize}
Assertion (a) is straightforward. Let us prove (b).

If $f(0)\ne 0$ then $\teps_k(\lambda'_x)(f){\equal}
\lambda'_x(\varepsilon_k\circ f){\equal}
\lambda(\eta_0\circ\varepsilon_k\circ f)
{\equal}\lambda(\varepsilon_{k-1}\circ\eta_0\circ f)
{\equal}(\teps_{k-1}(\lambda))'_{\teps_k(x)}(f)$.

If $f(0){\equal}0$ then $\teps_k(\lambda'_x)(f){\equal}
\lambda'_x(\varepsilon_k\circ f)
{\equal}\theta(\widetilde{\varepsilon_k\circ f}(x))
{\equal}\theta(\tf\circ\teps_k(x))
{\equal}((\teps_{k-1}\lambda)')_{\teps_k(x)}(f)$.
\pfend

\begin{rem}
In the case of an \et\ groupoid $G$ and cohomology groups
with coefficients in $G$-sheaves, Haefliger \cite{Hae79},
following Atiyah and Wall in the case of discrete groups
\cite{Ati-Wal67}, characterized the cohomology groups as the unique sequence
of functors $H^n$ such that
\begin{itemize}
\item[(a)] $H^0{\equal}\Gamma_{\text{inv}}$,
\item[(b)] $H^*$ admits
long exact sequences, and
\item[(c)] $H^n(G;I(\aA)){\equal}0$ for all $n\ge 1$ and for
each $G$-sheaf $I(\aA)$.
\end{itemize}
\end{rem}

Let us now prove
Proposition~\ref{prop:cech-sheaf}.
First, we note that $\aAb\hookrightarrow I(\aA)^\com$ canonically:
if $c\in \aA^n(U)$, then $\varphi(x){\equal}\teps_0^*[c(\teps_0(x))]
\in \aA^{n+1}_x$ is a section of $I(\aA)^n$ over $U$.

Using the uniqueness part in
Proposition~\ref{prop:derived-functors}, it suffices to show
that
$$\vH^n(\Mb;I(\aAb)){\equal}H^n(\Mb;I(\aA)^\com){\equal}\{0\}\quad \forall
n\ge 1.$$
This is true for $\vH^n$ thanks to
Lemma~\ref{lem:I(A)-acyclic}.

Define inductively a resolution
\begin{equation}\label{eqn:resolutionI}
0\to I(\aA)^\com\to \lL^{\com,0}\to
\lL^{\com,1}\to\cdots
\end{equation}
by $\lL^{\com,0}{\equal}I(I(\aA))^\com$ and
$\lL^{\com,q+1}{\equal}I(\lL^{\com,q}/\lL^{\com,q-1})$.
Since $\lL^{p,q}$ is flabby for all $p,q\ge 0$, the double
complex $K{\equal}(\lL^{p,q}(M_p))$ computes
$H^*(\Mb,I(\aA)^\com)$ (see above the introduction of this section).

The $E_2$-term with respect to the first filtration is
$E_2^{p,q}{\equal}H^pH^q(K)$. Since
$$0\to I(\aA)^p\to \lL^{p,0}\to \lL^{p,1}\to\cdots$$
is an exact sequence of flabby sheaves,
$$0\to \Gamma(M_p;I(\aA)^p)\to\Gamma(M_p;\lL^{p,0})\to\cdots$$
is exact, hence $E_2^{p,q}{\equal}0$ for $q\ge 1$ and
$E_2^{p,0}{\equal}H^p(\Gamma(M_\ast;I(\aA)^\ast))$.
Using again Lemma~\ref{lem:I(A)-acyclic} for the cover
$\uU_n{\equal}\{M_n\}$, we get $E_2^{p,0}{\equal}0$ for $p\ge 1$. Finally,
$H^n(\Mb;I(\aA)^\com){\equal}0$ for all $n\ge 1$.

\section{Invariance by Morita equivalence}
Let $G$ be a topological groupoid and $\aAb$ be an abelian sheaf on $G_\com$.
We will show that $H^*(G_\com;\aAb)$ and $\vH^*(G_\com;\aAb)$ are
invariant under Morita equivalence.

More precisely, if $G'$ is another groupoid and ${\aA'}^\com$ is a
sheaf on $G'_\com$, we say that $(G,\aAb)$ is Morita equivalent
to $(G',{\aA'}^\com)$ if there exists a groupoid $G''$, a sheaf ${\aA''}^\com$
on ${G''}^\com$ and (continuous) groupoid morphisms
$$G\stackrel{f}{\leftarrow}
G''\stackrel{f'}{\to}G'$$
such that $f$ and $f'$ are Morita equivalences and ${\aA''}^\com\cong
f^*\aAb\cong {f'}^*{\aA'}^\com$. Then

\begin{prop}
With the above assumptions, $f$ and $f'$ induce isomorphisms in
sheaf and \v{C}ech cohomology, thus
$$H^*(G_\com;\aAb)\cong H^*(G'_\com;{\aA'}^\bullet)\;\mbox { and }\;
\vH^*(G_\com;\aAb)\cong \vH^*(G'_\com;{\aA'}^\bullet).$$
\end{prop}

\pf
By standard arguments (compare with Proposition~\ref{prop:Morita}),
it suffices to show that for any open cover $\uU{\equal}(U_i)_{i\in I}$ of $G_0$,
the canonical morphism $f\colon G[\uU]\to G$ induces isomorphisms
$H^*(G_\com;\aAb)\cong H^*(G[\uU]_\com;f^*\aAb)$ and
$\vH^*(G_\com;\aAb)\cong \vH^*(G[\uU]_\com;f^*\aAb)$.
Below, we will abusively write $H^*(G[\uU]_\com;\aAb)$
instead of $H^*(G[\uU]_\com;f^*\aAb)$
\par\medskip
For \v{C}ech cohomology, using Remark~\ref{rem:cohomology-cover},
we have
\begin{equation}\label{eqn:cech1}
\vH^n(G_\com;\aAb){\equal}\lim_{\vV}\lim_{\wWb} H^n(\wWb;\aAb)
\end{equation}
where $\vV{\equal}(V_j)_{j\in J}$
runs over open covers of $G_0$ and $\wWb$ runs over open
covers of $G[\vV]$ such that $\wW_0$ consists of the single
open set $\coprod V_j$.

Similarly,
\begin{equation}\label{eqn:cech2}
\vH^n(G[\uU]_\com;\aAb){\equal}\lim_{\vV'}\lim_{\wW'_\com}
H^n(\wW'_\com;\aAb).
\end{equation}

where $\vV'{\equal}(V'_j)_{j\in J'}$
runs over open covers of $G[\uU]_0$ and $\wW'_\com$
runs over open covers of $G[\uU][\vV']$ such that $\wW'_0$
consists of the single open set $\coprod V'_j$.

Now, note that if $\vV'$ is an open cover
of $G[\uU]_0$ which is finer than the cover $(\{i\}\times U_i)_{i\in I}$,
then there exists an open cover $\vV$ of $G_0$ such that
$G[\vV]\cong G[\uU][\vV']$ (the elementary proof is left to the
reader). Therefore, in the right hand sides of (\ref{eqn:cech1}) and
(\ref{eqn:cech2}), the terms $\lim_{\wWb} H^n(\wWb;\aAb)$
and $\lim_{\wW'_\com}H^n(\wW'_\com;\aAb)$ are identical.
It follows that the right hand sides of (\ref{eqn:cech1}) and
(\ref{eqn:cech2}) are equal, hence \v{C}ech cohomology is invariant
by Morita equivalence.
\par\bigskip

From the above, we already find that sheaf cohomology is
invariant under Morita-equivalence when the groupoid is paracompact.
In fact, this holds for a general topological groupoid. Let us
sketch the proof for completeness.

Consider the resolution
$$\aAb\to \lL^{\com,0}\to \lL^{\com,1}\to\cdots$$
constructed like (\ref{eqn:resolutionI}).
Since $\lL^{p,q}$ is flabby for all $p$, $q$,
the double complex $(\lL^{p,q}(G_p))$ computes $H^*(G_\com;\aAb)$
and since the lines are exact (\ref{lem:I(A)-acyclic}),
$H^*(G_\com;\aAb)$ is the cohomology of the complex
$(\Gamma_{\mathrm{inv}}(\lL^{\com,q}))_{q\in \N}$.

Similarly, $H^*(G[\uU]_\com;\aAb)$ is the cohomology of the complex
$(\Gamma_{\mathrm{inv}}(f^*\lL^{\com,q}))_{q\in \N}$.
Now, it is elementary to check that
for every sheaf $\bBb$, $\Gamma_{\mathrm{inv}}(\bBb)$
is isomorphic to $\Gamma_{\mathrm{inv}}(f^*\bBb)$.
\pfend

\end{document}